\documentclass[11pt,leqno]{amsart}
\usepackage{amsmath,txfonts}
\usepackage[colorlinks=true,urlcolor=blue,
citecolor=red,linkcolor=blue,linktocpage,pdfpagelabels,
bookmarksnumbered,bookmarksopen]{hyperref}
\usepackage{epsfig,graphics,color}
\numberwithin{equation}{section}
\definecolor{verde}{rgb}{0,0.6,0}

\newcommand{\eps}{\varepsilon}
\newcommand{\vfi}{\varphi}

\newcommand{\la}{\lambda}

\newcommand{\GO}{{\mathcal G}_{\Omega}}
\newcommand{\G}{\mathcal G}
\newcommand{\MH}{\mathcal H}
\newcommand{\nh}{\mathcal N}
\newcommand{\nhlo}{{\mathcal N}_{\Omega}}

\newcommand{\io}{I_{\Omega}}

\newcommand{\hu}{\hat{u}}

\newcommand{\firt}{\Phi^{R\theta}}

\newcommand{\firo}{\Phi^{y}}

\newcommand{\firy}{\Phi^{Ry}}
\newcommand{\phiyn}{\Phi^{R_{n}y_{n}}}

\newcommand{\firxz}{\Phi^{Rx_{0}}}
\newcommand{\firxzn}{\Phi^{R_{n}x_{0}}}

\newcommand\R{\varmathbb{R}}

\newcommand\huno{H^{1}(\R^{N})}

\newcommand\homega{H^{1}_{0}(\Omega)}

\newcommand{\intr}{\int_{\R^{N}}}
\newcommand{\into}{\int_{\Omega}}

\newtheorem{theorem}{Theorem}[section]
\newtheorem{remark}[theorem]{Remark}
\newtheorem{lemma}[theorem]{Lemma}
\newtheorem{proposition}[theorem]{Proposition}
\newtheorem{corollary}[theorem]{Corollary}
\newtheorem{definition}[theorem]{Definition}
\newcommand{\bos}{\begin{remark}\rm}
\newcommand{\eos}{\end{remark}}
\newcommand{\bte}{\begin{theorem}}
\newcommand{\ete}{\end{theorem}}
\newcommand{\bpr}{\begin{proposition}}
\newcommand{\epr}{\end{proposition}}
\newcommand{\bdf}{\begin{definition}\rm}
\newcommand{\edf}{\end{definition}}
\newcommand{\bco}{\begin{corollary}}
\newcommand{\eco}{\end{corollary}}
\newcommand{\ble}{\begin{lemma}}
\newcommand{\ele}{\end{lemma}}
\newcommand{\bdm}{\begin{displaymath}}
\newcommand{\edm}{\end{displaymath}}
\newcommand{\beq}{\begin{equation}}
\newcommand{\eeq}{\end{equation}}
\newcommand{\bdim}{{\noindent{\bf Proof.}\quad}}
\newcommand{\edim}{{\unskip\nobreak\hfil\penalty50
\hskip2em\hbox{}\nobreak\hfil\mbox{\rule{1ex}{1ex} \qquad}
\parfillskip=0pt \finalhyphendemerits=0\par\medskip}}

\newcommand{\ov}{\overline}

\newcommand{\dys}{\displaystyle}
\begin{document}

\title[Solutions for Asymptotically Linear Problems in Exterior Domains]{\sc Positive Solutions for Asymptotically Linear Problems in Exterior Domains}

\author{Liliane A. Maia}
\address[L. A. Maia]{Departmento de Matem\'atica,
Universidade de Bras\'\i{}lia, 70910-900 Brasilia, Brazil.}
\email[L.A. Maia] { lilimaia@unb.br}
\thanks{Research of the first author is partially supported by CNPq/Brazil, PROEX/CAPES and FAPDF193.000.939/2015}
\author{Benedetta Pellacci}
\address[B. Pellacci]{Dipartimento di Scienze e Tecnologie, 
Universit\`a di Napoli ``Parthenope'', Centro Direzionale, Isola C4 80143 Napoli, Italy.}
\email[B. Pellacci]{benedetta.pellacci@uniparthenope.it}
\thanks{Research of the second author is partially 
supported by   MIUR-PRIN-2012-74FYK7 Grant: ``Variational and perturbative aspects of nonlinear differential problems'',
by the Italian Project FIRB 2012: ``Dispersive
dynamics: Fourier Analysis and Variational Methods", ``Gruppo
Nazionale per l'Analisi Matematica, la Probabilit\`a e le loro
Applicazioni'' (GNAMPA) of the Istituto Nazionale di Alta Matematica
(INdAM); University Project ``Sostegno alla ricerca individuale 
per il trienno 2015-2017".}

\begin{abstract}
\noindent
The existence of a positive solution for a class of 
asymptotically linear problems in exterior domains is 
established via a linking 
argument on the Nehari manifold and by means of 
a  barycenter function.
\end{abstract}
\maketitle


\section{Introduction}
In this paper we look for a positive solution for the following class
of elliptic problems
\begin{equation}\label{P}
\begin{cases}
-\Delta u + \lambda u= \displaystyle f(u), &
\mbox{in}\quad \Omega,
\\\\
\quad u\in H^{1}_{0}(\Omega)
\end{cases}
\end{equation}
where $\Omega$ is an exterior domain, i.e. a unbounded smooth domain 
in $\R^{N}$, $N\geq 2$, whose complement 
$\R^{N}\setminus\Omega$ is bounded 
and not empty.
The non-linearity  $f$ is super-linear at zero and
asymptotically linear at infinity, and
the model example we have in mind is 
\beq\label{fmod}
f(t)=\frac{t^{3}}{1+st^{2}},\qquad 0<s<\frac1\la.
\eeq
This kind of equations arises in the study of the propagation of a 
light  beam in saturable media, as photorefractive crystals,
which are largely exploited
in experimental observations because they require 
low optical power to exhibit very strong non-linear effects
(\cite{ak98,lit,oski,stse}). 
When the medium has some additional properties (\cite{ak99}), 
the non-linearity given in \eqref{fmod} is approximated by
putting $s=0$ and obtaining  $f(t)=t^{3}$, that is the 
more treated non-linearity, also associated  to the study of the 
propagation of a light beam 
in the so-called Kerr media (\cite{ak98}).
From a mathematical point of view, this kind of problems 
has been extensively studied in the last thirty years, 
especially when $f(t)=|t|^{p-2}t$, and $p>2$ is subcritical.
The main difficulty relies in the lack of compactness that
can be overcome if $\Omega$ is radially symmetric and one looks
for radially symmetric solutions (\cite{beli, esli}). 
While, if $\Omega$ has no symmetry properties, the lack of
compactness becomes prevalent, moreover, it has been
observed in \cite{bece} that there are no least action solutions, 
so that  a higher action level solution 
has to be searched. As a consequence, a 
careful analysis of a general Palais-Smale sequence,
even not minimising, is needed and this has been done
in \cite{bece} introducing the nowadays 
so-called  ``splitting lemma''.
This analysis permits to locate the action levels where the Palais-Smale
condition holds, so that, the difficulty becomes to construct
a minimax level where this compactness property is
satisfied. In \cite{bece} 
suitable subspaces  of the $L^{p}(\Omega)$ sphere 
have been considered, with additional  prescriptions, in 
order to increase the minimax level. This approach led to the 
existence of a constrained minimisation point, which, thanks to
the homogeneity property of $f(t)=|t|^{p-2}t$ gives 
a positive solution when $\R^{N}\setminus \Omega$ has 
small diameter. This kind of argument has been developed in 
several directions: in \cite{balio}
the assumption on the size of $\R^{N}\setminus \Omega$ has been 
removed, in \cite{cepa92, cepa95}   multiplicity results have 
been obtained, in dependence on $\la$ and
on  the topology of $\Omega$ (see also the interesting 
survey \cite{ce});
in \cite{cecl, clsa}  sign-changing 
solutions have been found. 
\\
When $f$ has not a polynomial dependence with respect to $t$,
there are only a few contributions in the literature. In particular,
in the case in which $f$ is asymptotically linear,
in \cite{ci} the existence of a solution is proved
for $f$ depending on the spatial variable 
and satisfying some  assumptions corresponding, 
in the autonomous case,  to a convexity property. 
More recently, Problem \eqref{P} has been studied in
\cite{lizh} where it is found a positive solution
under some suitable hypotheses on $f$, 
among which, again the convexity for $t\in(0,+\infty)$
and  in \cite{mamiso}, where it is proved the existence 
of a  radial positive solution in the exterior of a ball.
\\
Notice that, $f$ given in \eqref{fmod}
is not convex so that the result in 
\cite{ci, lizh} do not cover the physical meaningful case.
Moreover, we do not suppose any symmetry property on 
the domain, so that there is no reason to expect a 
radially symmetric  solution. With this respect, our result 
in the model case \eqref{fmod} is related to the interest in  
the physical context of  asymmetric solutions that are observed
in the interaction between  the crystal and the spatial solitons 
or in the collision of  symmetric solitons (\cite{ak98, lit}). 
\\
Our existence result in  the model case
can be stated as follows.
\bte\label{rhointro}
Assume $N\geq 2$, and that $f$ is given in \eqref{fmod}. 
Then, Problem \eqref{P}
has at least a positive solution.
\ete
This theorem will follow as a consequence of
a more general  Theorem \ref{rho},
proved for a  general non-linearity $f(t)$,
(see Section \ref{setting}). 
\\
As in the super-linear case,  when $f(t)$ is given in \eqref{fmod},
one expects that
the least action level of the associated functional 
$\io:H^{1}_{0}(\Omega)\mapsto \R$
\beq\label{defio}
I_{\Omega}(u):=\into \dfrac12\left[|\nabla u|^{2}
+\la u^{2}\right]-F(u),\quad \text{with }\; F(u)=\int_{0}^{u}f(t)dt
\eeq
is equal to the least action level of the functional associated to the problem in the whole $\R^{N}$
$$
I(u):=\intr \dfrac12\left[|\nabla u|^{2}+\la u^{2}\right]-F(u),
$$
which is showed to be attained by a unique, positive,
radially symmetric  solution
(\cite{beli,begaka,seta}).
But, when dealing with asymptotically 
linear $f(t)$, even in the model case \eqref{fmod}, one 
loses the precious homogeneity property, crucial in order 
to reabsorb the Lagrange multiplier when finding
the least action solution as a constrained minimum point 
on a $L^{p}(\Omega)$ sphere. So that, working 
on $L^{p}(\Omega)$ spheres is not suitable in this context. 
In the light of these considerations, we will work on 
the well-known Nehari manifolds 
associated to the functionals $\io$ and $I$
(\cite{ne1,ne2})
\beq\label{nehari}
\begin{split}
\nhlo &=\{u\in H^{1}_{0}(\Omega)\setminus\{0\}:
\langle \io'(u),u\rangle=0\},
\\
{\mathcal N}&= \{ u \in H^1(\mathbb{R}^N)\setminus\{0\}
: 
\langle I'(u),u\rangle =0\},
\end{split}
\eeq
realizing the least action levels as minimum values of $\io$ and
$I$ on these manifolds
\begin{equation}\label{m.omega}
m_{\Omega}:=
\inf_{u\in {\mathcal{N}_{\Omega}}} I_{\Omega}(u),\qquad
m:=\inf_{u\in \mathcal{N}} I(u).
\end{equation}
This approach yields some new difficulties and some new advantages;
concerning the formers, notice that not every function in $H^{1}_{0}(\Omega)$ 
can be projected on $\nhlo$, so that, one has to be careful 
when  defining the projection of a function on $\nhlo$; 
But, thanks to a
monotonicity condition of $f$ (\eqref{fmono})
the projection turns out to be  unique when it exists. 
With this approach we will  show that 
$m_{\Omega}=m$ so that a higher energy critical level 
has to be searched also in this case. 
On the other hand, we have some benefits when looking for 
a compactness property; indeed, recall that, usually when
dealing with  asymptotically linear non-linearity the 
suitable notion of compactness property is the Cerami condition
introduced in \cite{ce2}; here, strongly exploiting the features
of $\nhlo$, we will be able to show that the usual Palais-Smale 
condition on $\nhlo$ holds.
However,  in doing this, we will again  have to handle the
lack of homogeneity of $f$. 
\\
Finally,  we will be able to choose 
linking sets as proper subsets of 
$\mathcal N_{\Omega}$.
At this stage, it will be important to show proper estimates
of $\io$ on these subsets. This  will be done by introducing 
a suitable abstract asymptotic threshold  
(in Lemma \ref{lemma:eps}), that will be  crucial 
in order to compare all the terms in the functional $\io$
(defined in \eqref{defio}).
This threshold generalises  the commonly used one in the context of
polynomial non-linearity.
Concluding, let us observe that, in order to highlight the main 
novelties and difficulties one encounters when facing this kind
of problems, we have preferred to study the existence of a solution,
but we believe that with the same tools introduced here
it  could be possible to face the question of multiplicity of solutions 
as done in \cite{cemopa, cepa95} for the 
polynomial case. Moreover, our argument can be used  to
prove an analogous result for a general,
not necessarily homogeneous, super-linear $f$
(see for more details Remark \ref{super}).
\\
Concluding this introduction, let us point out that the
use of the constraint of Nehari manifold has extensively
been used  in the related context of the problem in the whole 
space with the
presence of a  potential (see for instance 
\cite{amma,wi, szwe} and the references therein). 
In particular, in \cite{szwe,evwe} 
the properties of a generalized Nehari manifold
are deeply exploited in order to study indefinite problems;
here we choose to assume all the regularity conditions
needed to use the classical Nehari manifold, but 
it would be  interesting to use the analysis of \cite{szwe,evwe} 
to weaken our regularity assumptions on the non-linearity.
\\
The paper is organized as follows: in Section \ref{setting}
we settle our general context and state our main existence
result (Theorem \ref{rho}); 
in Section \ref{asymptotic} we prove some useful technical
propositions, while in Section \ref{compactness} we 
deal with the compactness property.  
Finally, in Section \ref{final} we conclude the proofs of 
Theorem \ref{rho} and \ref{rhointro}.
\section{Setting of the problem and main results}\label{setting}
We will work in $H^{1}_{0}(\Omega)$ and in
$H^{1}(\R^{N})$. The norms on these spaces will
be denoted respectively by 
$$
\|u\|^{2}_{\Omega}=(u,u)_{\Omega}=
\|\nabla u\|^{2}_{2,\Omega}+\la\|u\|^{2}_{2,\Omega},\quad
\|u\|^{2}=(u,u)=\|\nabla u\|^{2}_{2}+\la\|u\|^{2}_{2},
$$ where
$\|\cdot\|_{p.\Omega}$
($\|\cdot\|_{p}$) is the usual norm in $L^{p}(\Omega)$
(in $L^{p}(\R^{N})$).
\\
Every solution of \eqref{P}, 
is  a critical point of the $C^{1}$ functional $I_\Omega$,  
defined  in \eqref{defio}.
We have introduced the model example in \eqref{fmod}, 
however, in general, we will assume that  $f$ has the following properties
\beq
\label{fprol} 
\text{$f\in C^{1}[0,+\infty)\cap C^{3}(0,+\infty),$ } 
\text{ $f(t)\equiv 0,$ }  \text{ $\forall\, t\in (-\infty,0)$};
\eeq
there exists a positive constant $D$ and $p_{1}\leq p_{2}\in 
(1,\bar{p})$ where $\bar{p}=+\infty$ if $N=2$ and 
$\bar{p}= (N+2)/N-2$ if $N\geq 3$,  such that
\beq\label{fdergrowth}
|f^{(k)}(t)|\leq D\left[|t|^{p_{1}-k}+|t|^{p_{2}-k}\right],\;k=0,\,1,\,2,\,3.
\eeq
The following assumption will be crucial in all our arguments:
\beq\label{fmono} 
\frac{f(t)}{t} \;\text{is  an} 
\text{ increasing function}  \; \forall\,t >0.
\eeq
Moreover, we suppose that  $f$ is asymptotically linear at infinity,
namely, we assume that there exists
a positive number $l_{\infty}$ such that 
\beq\label{finfty}
\displaystyle \lim_{t\rightarrow +\infty}\frac{f(t)}{t} =l_{\infty} \,,
\qquad \lambda- l_{\infty} < 0\,,
\eeq
where the second inequality is assumed 
in order to have a solution of the problem in the whole $\R^{N}$
(see \cite{beli, begaka, stzh}).
When dealing with this kind of non-linearities
it is often assumed a non-quadraticity type condition
(introduced in \cite{coma})
\beq\label{NQ}
f(t)t - 2F(t) \geq 0, \text{ for all }  t \in \mathbb{R^{+}}
\quad \text{and}\quad 
\displaystyle \lim_{t \to +\infty}[f(t)t - 2F(t)]= +\infty\; .
\eeq
Finally, in order to let the Nehari manifold be a natural constraint, 
we will assume that
\beq\label{iponh}
f'(t)-\frac{f(t)}t> 0\quad \text{for every  $t>0$}.
\eeq
\begin{remark}
Notice that all the previous hypotheses are satisfied 
in the model case \eqref{fmod}.
Moreover note that,
as a consequence of \eqref{fmono} and \eqref{finfty}, 
the general non-linearity $f(t)$ is always below the line 
$l_{\infty} t$, for every $t>0$.
\end{remark}
In order to prove our existence result we will make a comparison
with the following problem in the whole $\R^{N}$ which will be
called problem at infinity
\begin{equation}\label{Pinf}
\begin{cases}
-\Delta u + \lambda u= \displaystyle f(u), &
\mbox{in}\quad \R^{N},
\\\\
\displaystyle \lim_{|x|\rightarrow \infty}u(x)=0.
\end{cases}
\end{equation}
The problem \eqref{Pinf} has a positive, radially symmetric,
least action solution (whose existence is proved in \cite{beli,stzh}  
for $N\geq 3$ and  in \cite{begaka} for $N= 2$.), which we denote with $w$. In \cite{seta} it is shown the uniqueness of $w$ when
$f$ satisfies the additional hypothesis
\[
g(t):=\frac{-\la t+tf'(t)}{-\la t+f(t)}\,\;\text{ is decreasing in 
$[b,+\infty)$ }
\]
where $b$ is the (unique thanks to \eqref{fmono}) 
positive number given by $f(b)=\la b$. 
Notice that this condition is satisfied in the case of \eqref{fmod}.
Since  we will also  deal with a more general non-linearity and 
will not  assume this hypothesis we will suppose that
\beq\tag{$U$}\label{unique}
\text{The  positive radially symmetric solution of  problem (\ref{Pinf}) is unique.} 
\eeq
Our main result is the following.
\bte\label{rho}
Assume $N\geq 2$ and hypotheses  \eqref{fprol}, \eqref{fdergrowth},
\eqref{fmono}, \eqref{finfty}, \eqref{NQ}, \eqref{iponh} and 
\eqref{unique}. 
Then Problem \eqref{P} has at least a positive solution.
\ete
\bos
In \cite{ci} Problem \eqref{P} is studied for 
a non-autonomous $f$ satisfying some 
suitable hypotheses including one, that becomes, in
the autonomous case, a convexity condition.
In \cite{lizh} Problem \eqref{P} is studied assuming \eqref{fprol},
\eqref{fmono}, \eqref{finfty}, \eqref{NQ}, \eqref{iponh}, 
\eqref{unique} plus other
conditions among which it is assumed that
f  is a convex function. 
Note that,  $f$ given in \eqref{fmod}
 is not convex, so that the existence
of a positive solution in exterior domains for this 
kind of non-linearity cannot be deduced from the result
of \cite{ci, lizh}.
Obviously, our result applies also to non-linearities
that have infinitely many flex points.
\eos
\bos\label{super}
Our argument can be exploited to deal with super-linear, not necessarily
homogeneous, non-linearities assuming conditions \eqref{fprol},
\eqref{fdergrowth}, \eqref{fmono}, \eqref{iponh} and, in place
of \eqref{finfty} and \eqref{NQ}, supposing the classical
Ambrosetti-Rabinowitz condition
\[
\exists \, \mu>2,\, :\,0<\mu F(t)\leq tf(t),\quad \forall\,t\in \R^{+}\setminus \{0\}.
\]
Indeed, also under these hypotheses, it is possible to exploit
the properties of  the Nehari manifold $\nhlo$
 as explained in \cite{ra}.
\eos
As already said in the introduction, we will work on the Nehari
manifold introduced in \eqref{nehari}. Nowadays, this has become
a classical tool in variational methods thanks of  
its useful features as it has been also highlighted in the recent
contribution \cite{nove}.
The following remark clarifies the role of \eqref{iponh} in order 
to use $\nhlo$.

\bos \label{natu}
Let us first observe that $\nhlo$ is 
the set of non-trivial zeroes of the function 
$N_{\Omega}:H^{1}_{0}(\Omega)\setminus\{0\}\mapsto \R$
given by
\[
N_{\Omega}(u)=\langle I'_{\Omega}(u),u\rangle=\|u\|_{\Omega}^{2}-\into f(u)u.
\]
Notice that $\nhlo$ is
actually a manifold and it is a natural constraint. 
Namely,  for every $u\in\nhlo$ it results
\(
\langle N'_{\Omega}(u),u\rangle<0.
\)
Indeed, consider  $u\in \nhlo$  and 
use \eqref{fdergrowth} and \eqref{iponh} to obtain
\begin{align*}
\langle N_{\Omega}'(u),u\rangle
&=
2\|u\|^{2}-\into [f'(u)u^{2}+f(u)u]=
\into u^{2}\left[\frac{f(u)}{u}-f'(u)\right]< 0.
\end{align*}
Now, suppose that $u\in\nhlo$ is a constrained critical point
of $\io$, then there exists a real number $\mu$ such that
$\io'(u)-\mu N_{\Omega}'(u)=0$; taking $u$ as test function 
one gets
$\mu \langle N_{\Omega}'(u),u\rangle=0$
then yields $\mu=0$, i.e. $u$ is a free critical point.
\eos

\section{Asymptotic Estimates}\label{asymptotic}
In the sequel, when we compute the $H^{1}(\R^{N})$ norm of a function 
$u\in H^{1}_{0}(\Omega)$, it is implicitly assumed that
$u$  is extended as zero outside $\Omega$ so that $u\in H^{1}(\R^{N})$. 
Let us introduce a cut-off $C^\infty$ function $\xi: \mathbb{R}^N \to [0,1]$  defined by
\beq
\xi (x):= \tilde \xi \displaystyle{\left(\frac{|x|}{\rho}\right)}\;,
\eeq
with $\rho$  being the smallest positive number such that $\mathbb{R}^N \backslash \Omega \subset B_\rho(0)$ 
and $\tilde \xi : \mathbb{R}^+ \cup \{0\} \to [0,1]$ being a non decreasing function such that
\beq\label{defxi}
\tilde \xi(t)=0 \quad \text{if} \quad t \leq 1 \quad \text{and} \quad \tilde \xi (t)= 1 \quad \text{if} \quad t \geq 2,\quad 
| \tilde{\xi}'(t)| \leq 2.
\eeq
Recall that $w$, the least action positive solution of 
Problem \eqref{Pinf}  enjoys the following
asymptotic behavior (see \cite{beli})
\beq\label{behavior:w}
\lim_{R\to +\infty}w(R)R^{(N-1)/2}e^{\sqrt{\lambda}R}=\sigma,
\qquad
\lim_{R\to +\infty}w'(R)R^{(N-1)/2}e^{\sqrt{\lambda}R}=-\sigma\sqrt{\lambda},
\eeq
we will often use the following Lemma proved in Lemma II.2 
in \cite{balio},
(see also Proposition 1.2 in \cite{bali}).
\begin{lemma}\label{bali}
Let $\vfi_{1}\in C(\R^{N})\cap L^{\infty}(\R^{N}),\,\vfi_{2}\in C(\R^{N})$ satisfy for some
$\alpha,\,\beta \geq 0$, $c\in \R$
$$
\lim_{|x|\to+\infty}\vfi_{1}(x)e^{\alpha|x|}|x|^{\beta}=c,
\qquad
\intr|\vfi_{2}(x)|e^{\alpha|x|}(1+|x|^{\beta})dx<+\infty,
$$
then
$$
\lim_{|\zeta|\to +\infty}e^{\alpha|\zeta|}|\zeta|^{\beta}\intr  \vfi_{1}(x+\zeta)\vfi_{2}(x)dx=c\intr \vfi_{2}(x)dx.
$$
\end{lemma}
In what follows we will  use the notation
\beq\label{def:tras}
w^{\theta}(\cdot):=w(\cdot-\theta),\qquad \forall\,\theta\in \R^{N}\setminus \{0\}.
\eeq
As a consequence of the previous lemma,  it is easy to prove the following result.
\begin{lemma}\label{lemma:eps}
Assume \eqref{fdergrowth} and let 
$x_{0}\in \R^{N}$ with $|x_{0}|=1$ and $y\in \partial B_{2}(x_{0})$.
Let us define the quantity  
\[
\eps_{R}:=\intr f(w^{Rx_{0}})w^{Ry}=\intr 
f(w^{Ry})w^{Rx_{0}}.
\]
Then, $\eps_{R}$ satisfies 
\[
\lim_{R\to+\infty}\eps_{R}
\left[ (2R)^{(N-1)/2}e^{ 2 R\sqrt{\la}}\right]=c_{0}>0.
\]
\end{lemma}
\begin{remark}
Let us observe that, in 
\cite{evwe}  an analogous asymptotic  threshold  
has been introduced in order to deduce precise energy estimates,
for an indefinite problem in the whole space with 
super-linear, but not necessarily homogeneous, non-linearities.\end{remark}
\bdim
The proof can be done following \cite{bali} or \cite{balio}; indeed by a 
change of variable, it results
$$
\eps_{R}=\intr f(w)w^{R(x_{0}-y)}.
$$
In order to obtain the conclusion it is sufficient to apply 
Lemma \ref{bali} with 
$\vfi_{1}=w,\,\vfi_{2}=f(w)$ and $\zeta=-R(x_{0}-y)$.
Let us check that all the hypotheses are satisfied.
Using  \eqref{behavior:w}, we obtain that $\vfi_{1}$ satisfies
the hypothesis of Lemma \ref{bali}
with $\alpha =\sqrt{\la}$ and $\beta=(N-1)/2$.
Let us see what happens with $\vfi_{2}$. 
The first limit in \eqref{behavior:w} and condition \eqref{fdergrowth} 
imply that there exists $R_{1}$ such that
$$
\vfi_{2}(x)=f(w)\leq 2D\sigma\left[
e^{-\alpha p_{1}|x|}|x|^{-\beta p_{1}}+e^{-\alpha p_{2}|x|}
|x|^{-\beta p_{2}}\right]\quad \forall\, |x|>R_{1}.
$$
Then it holds, with $C$ a positive constant
\begin{align*}
\intr f(w(x))e^{\alpha |x|}(1+|x|)^{\beta}
&\leq C\text{meas}B_{ R_{1}}(0)
+DC\int_{\{|x|>R_{1}\}}
e^{(1- p_{1})\alpha|x|}|x|^{-\beta p_{1}}
(1+|x|^{\beta})
\\
&+DC\int_{\{|x|>R_{1}\}}
e^{(1- p_{2})\alpha|x|}|x|^{-\beta p_{2}}(1+|x|^{\beta})
\end{align*}
so that the hypotheses of Lemma \ref{bali}
are satisfied since $p_{i}>1$.
\edim
Moreover, for every $q>1$, for every 
$\theta\in \R^{N}\setminus\{0\}$ and for every compact set $K$
Lemma \ref{bali} implies
\beq\label{opic}
\int_{K}\left[w^{R\theta}\right]^{q}
\leq O\left(\frac1{R^{(N-1)q/2}e^{Rq\sqrt{\la}}}
\right).
\eeq
Indeed,  denote with $a$ the modulus of $\theta$ and observe 
that \eqref{behavior:w} implies that there exists $\rho_{0}$ such that
 it results
\beq\label{modz}
w(z)(|z|)^{(N-1)/2}e^{\sqrt{\la}|z|}\leq 2\sigma,
\quad \text{for every $|z|\geq \rho_{0}$}.
\eeq
Let us fix $r_{0}$ such that $K\subseteq B_{r_{0}}(0)$
and put $z=x-R\theta$.
We have
\[
\begin{split}
R^{(N-1)q/2}e^{Rq\sqrt{\la}}\int_{K}\left[w^{R\theta}\right]^{q}
\leq
R^{(N-1)q/2}e^{Rq\sqrt{\la}}\int_{B_{r_{0}}(R\theta)}\left[
w(z)\right]^{q}dz
\end{split} 
\]
Moreover, $|z|\geq Ra-r_{0}$, so that,
for every $R$ such that $Ra-r_{0}>\rho_{0}$,   
inequality \eqref{modz} implies that there exists
a positive constant $C$ such that
\[
\begin{split}
R^{(N-1)q/2}e^{Rq\sqrt{\la}}\int_{K}\left[w^{R\theta}\right]^{q}
\leq
(2\sigma)^{q}C.
\end{split} 
\]
Therefore, for every $q\geq 2$, for every 
$\theta\in \R^{N}\setminus\{0\}$ and for every compact set $K$,
it results
\beq\label{compact}
\int_{K}\left[w^{R\theta}\right]^{q}
\leq o(\eps_{R}),\qquad 
\int_{K}\left[\left|\nabla w^{R\theta}\right|\right]^{q}
\leq o(\eps_{R}), 
\eeq
where the second inequality follows by an analogous 
argument for $w'$.
For every $\theta\in \R^{N}\setminus\{0\}$  
let us define the map $ \Phi^{\theta} : \mathbb{R}^{+} \to H^1_0(\Omega)$ by
\begin{equation} \label{phi}
R\mapsto \Phi^{\theta}(R)(\cdot):= \xi(\cdot)  w(\cdot-R\theta).
\end{equation}
In the sequel, for simplicity, we will use the notation
\[
\firt:=\Phi^{\theta}(R)
\]
\begin{lemma}\label{phiprop}
Let us assume \eqref{fprol} and  \eqref{fdergrowth}.
Let $\theta\in \R^{N}\setminus \{0\}$.
Then,  the function $\firt$ is continuous in $R$  and 
\begin{align}\label{limy}
\|\firt-w^{R\theta}\|^{2} \leq o(\eps_{R}), \; &
\quad
\left|\|\firt\|^{2}-\|w^{R\theta}\|^{2}\right| \leq o(\eps_{R}), \;
\\
\label{limyI}
\left| I_\Omega(\firt) - I(w)\right|&\leq o(\eps_{R}).\; 
\end{align}
Moreover, for every $\tau \geq 0$ it results
\beq\label{Fprop}
\left|\intr F(\tau \firt) -F(\tau w^{R\theta})\right|\leq 
\left[\tau^{p_{1}+1}+\tau^{p_{2}+1}\right] o(\eps_{R}).
\eeq
\end{lemma}
\bdim
The continuity with respect to $R$ of $\firt$ 
follows from the regularity properties 
of $w$. Moreover, taking into account \eqref{defxi},
\eqref{def:tras} and using Young inequality, it follows 
\begin{align*}
\|\firt-w^{R\theta}\|^{2}=& 
\intr |\nabla \xi w^{R\theta}+(\xi-1)\nabla w^{R\theta}|^{2}+
\la \intr |\xi-1|^{2}(w^{R\theta})^{2}
\\
\leq&  \int_{B_{2\rho}(0)}
(w^{R\theta})^{2}
\left[ |\nabla \xi|^{2}|+\la  |\xi-1|^{2}\right]
+\int_{B_{2\rho}(0)}
 |\xi-1|^{2}|\nabla w^{R\theta}|^{2}
\\
&+2\int_{B_{2\rho}(0)}w^{R\theta}|\xi-1||\nabla \xi|
\left|\nabla w^{R\theta}\right|
\\
\leq&  \int_{B_{2\rho}(0)} 
(w^{R\theta})^{2}
\left[ 2|\nabla \xi|^{2}|+\la  |\xi-1|^{2}\right]
+2\int_{B_{2\rho}(0)}
 |\xi-1|^{2}|\nabla w^{R\theta}|^{2}
\\
\leq& 
2\int_{B_{2\rho}(0)} 
|\nabla w(x-R\theta)|^{2}dx+
C\int_{B_{2\rho}(0)} |w(x-R\theta)|^{2}dx,
\end{align*}
with $C$ a positive constant.
Then, \eqref{limy} follows from \eqref{compact}.
\\
In order to show the second conclusion in \eqref{limy}, note that \eqref{defxi} 
and \eqref{compact} imply
\begin{align*}
\left|\left(w^{R\theta}, \firt-w^{R\theta}\right) \right|
&\leq 
\int_{B_{2\rho}(0)}\!\!3\left|\nabla w^{R\theta}\right|^{2}+
2( w^{R\theta})^{2}\leq o(\eps_{R});
\end{align*}
this and the first conclusion in \eqref{limy} yield the second one.
\\
To prove \eqref{limyI}, let us observe that
$I_{\Omega}(\firo)=I(\firo)$. Therefore,
from \eqref{defxi}, \eqref{fdergrowth}, 
and using Lagrange mean value 
Theorem, we get that there exists a positive 
constant $C_{0}$ such that
\begin{align*}
|I_\Omega(\firt) - I(w)|&\leq 
\left| \|\firt\|^{2}-\|w^{R\theta}\|^{2}\right|
+\left|\intr F(\xi w^{R\theta})-F(w^{R\theta}) 
\right|
\\
&\leq  \left| \|\firt\|^{2}-\|w^{R\theta}\|^{2}\right|
+C_{0}D\int_{B_{2\rho}(0)} (w^{R\theta})^{p_{1}+1}+
(w^{R\theta})^{p_{2}+1}.
\end{align*}
Then, \eqref{limyI} follows from
\eqref{compact}  and \eqref{limy}.
Moreover, 
taking into account \eqref{defxi} and \eqref{def:tras},
and using condition \eqref{fdergrowth} 
we get that there exists a positive constant $C_{1}$ such that
\begin{align*}
\left|\intr  F(\tau w^{R\theta})-F(\tau \firt)\right|
&\leq 
C_{1}D
\int_{B_{2\rho}(0)}\tau^{p_{1}+1}(w^{R\theta})^{p_{1}+1}
+\tau^{p_{2}+1}(w^{R\theta})^{p_{2}+1}.
\end{align*}
Then the conclusion follows from \eqref{compact}.
\edim
As already said, the linking sets will be subsets of $\mathcal N_{\Omega}$, but, as $f$ is asymptotically linear, not
every function of $H^{1}_{0}(\Omega)$ can be projected 
on $\mathcal N_{\Omega}$.
The following lemma tells us when this projection can be
performed.

\begin{lemma}\label{tprop}
Let us assume \eqref{fprol}, \eqref{fmono}, \eqref{finfty},  
and define the maps
$$
\GO(u)=\|\nabla u\|_{2,\Omega}^{2}+\lambda\|u\|_{2,\Omega}^{2}
-l_{\infty}\|u^{+}\|_{2,\Omega}^{2}\,,\qquad
\G(u)=\|\nabla u\|_{2}^{2}+\lambda \|u\|_{2}^{2}-l_{\infty}
\|u^{+}\|_{2}^{2}\,.
$$
Then, the maps 
\begin{align*}
T_{\Omega}:A_{\Omega}&:=\left\{u \in H^{1}_{0}(\Omega)\,:\,\GO(u)<0\right\}\mapsto \R^{+},\quad
u\mapsto T_{\Omega}(u)\,:\,
T_{\Omega}(u)u\in \nhlo ,
\\
T\,:A &:=\left\{u \in H^{1}(\R^{N})\,:\,\G(u)<0\right\}
\mapsto \R^{+},\quad\;\;
u\mapsto T(u)\;\,:\,
T(u)u\in \nh
\end{align*}
are well defined and continuous.
In addition, for every  $\theta\in\R^{N}\setminus\{0\}$,
there exists $R_{0}=R_{0}(|\theta|)$  such that 
$\firt\in A_{\Omega}\cap A$  for every $R\geq R_{0}$ and 
\beq\label{limtn}
\lim_{R\to +\infty}T_{\Omega}(\firt)= 1.
\eeq
\end{lemma}
\bdim
Let us define $g_{\Omega}\,: [0,+\infty)\times A_{\Omega}\mapsto \R$ by
$$
g_{\Omega}(\tau,u):=\begin{cases}
\| u\|_{\Omega}^{2}=\|\nabla u\|_{2,\Omega}^{2}+\la \|u\|^2_{2,\Omega}& \tau=0,
\\\\
\dys\frac1{\tau^{2}}\langle I'_{\Omega}(\tau u),\tau u\rangle & \tau>0.
\end{cases}
$$
Observe that $g_{\Omega}$ is continuous on $\R^{+}\times A_{\Omega}$, 
moreover, hypotheses \eqref{fprol}, \eqref{fmono} and  
\eqref{finfty} allow to apply
Lebesgue Dominate Convergence Theorem to 
obtain
$$
\lim_{\tau \to 0^{+}} g_{\Omega}(\tau,u)= \| u\|_{\Omega}^{2}-
\lim_{\tau\to 0^{+}} \into\frac{f(\tau u)\tau u}{\tau^{2}}=
\| u\|_{\Omega}^{2}-\lim_{\tau\to 0^{+}} \into\frac{f(\tau u)}{\tau u}u^{2}=\| u\|_{\Omega}^{2}\,.
$$
So that $g_{\Omega}$ is continuous on $[0,+\infty)\times A_{\Omega}$
and $g_{\Omega}(0,u)>0$. While, \eqref{fprol}, \eqref{fmono}
and \eqref{finfty}  yield
$$
\lim_{\tau \to +\infty} g_{\Omega}(\tau,u)=
\| u\|_{\Omega}^{2}-
\lim_{\tau\to +\infty} \into\frac{f(\tau u)}{\tau u}u^{2}=
\|u\|_{\Omega}^{2}
-l_{\infty}\|u^{+}\|_{2,\Omega}^{2}<0\,,
$$
where the last inequality is implied by the inequality $\GO(u)< 0$.
Finally, \eqref{fmono} implies that $g_{\Omega}$ is strictly
decreasing with respect to $\tau$, so that,
since it is a continuous function, there exists a unique
$T_{\Omega}(u) > 0$ such that $g_{\Omega}( T_{\Omega}(u),u) =0$, that is
\[
\langle I_{\Omega}^\prime(T_{\Omega}(u) u), T_{\Omega}(u) u \rangle =0,
\]
i.e. $T_{\Omega}(u) u \in \mathcal{N}_{\Omega}$, showing 
that $T_{\Omega}$ is well defined. \\
In order to show that the map $T_{\Omega}$ is continuous 
in $A_{\Omega}$, let $u$ be in $A_{\Omega}$, so that 
we can consider $\tilde{\Omega}\subset A$ with $|\tilde{\Omega}|>0$ and $u^{+}(x)\neq 0$ iff  $x\in \tilde{\Omega}$ and
assume that $(u_{n})\in A_{\Omega}$
converges to  $u$.
The continuity property of $\GO$ implies that 
there exists $n_{0}$ such that $\GO(u_{n})<0$ for every $n\geq 
n_{0}$, so that there exists $T_{\Omega}(u_{n})$ which will
be denoted as $T_{n}$.
Assume, by contradiction that, up to a subsequence, $T_{n}\to +\infty$, then 
\eqref{fprol}, \eqref{fmono} and \eqref{finfty} 
allow to apply Lebesgue dominated convergence Theorem to obtain 
that
$$
\lim_{n\to +\infty}\int_{\tilde{\Omega}}\frac{f(T_{n}u_{n})}{T_{n}u_{n}}u^{2}_{n}=
l_{\infty}\int_{\tilde{\Omega}}u^{2}=l_{\infty}\into (u^{+})^{2}.
$$
Then, by definition of $T_{n}=T_{\Omega}(u_{n})$, 
and using \eqref{fprol}, \eqref{fmono} and \eqref{finfty} 
it results
\begin{align*}
\|\nabla u\|^2_{2,\Omega}+\la \|u\|^2_{2,\Omega}&=\lim_{n\to \infty}
\|\nabla u_{n}\|^2_{2,\Omega}+\la \|u_{n}\|^2_{2,\Omega}= 
\lim_{n\to \infty}\into\frac{f(T_{n}u_{n})}{T_{n}u_{n}}u_{n}^{2}
\\
&=\lim_{n\to \infty}\int_{\tilde{\Omega}}\frac{f(T_{n}u_{n})}{T_{n}u_{n}}u_{n}^{2}=l_{\infty}\int_{\tilde{\Omega}}u^{2}
= l_{\infty}\into (u^{+})^{2}
\end{align*}
implying  $\GO(u)=0$, which is a contradiction. Then, up to a subsequence,
$T_{n}\to T_{0}$; if it were $T_{0}=0$ then \eqref{fprol},
\eqref{fmono} and \eqref{finfty}  would imply
\begin{align}\label{tconv}
\|\nabla u\|^2_{2,\Omega}+\la \|u\|^2_{2,\Omega}&=\lim_{n\to \infty}
\|\nabla u_{n}\|^2_{2,\Omega}+\la \|u_{n}\|^2_{2,\Omega}=\lim_{n\to \infty} 
\into\frac{f(T_{n}u_{n})}{T_{n}u_{n}}u_{n}^{2}=0,
\end{align}
again a contradiction, as $\GO(u)<0$.Then, $T_{n}\to T_{0}>0$
and  passing to the limit in \eqref{tconv}
we derive that $T_{0}=T_{\Omega}(u)$, as $T_{\Omega}(u)$ is 
unique.
In addition, take $\theta\in \R^{N}\setminus\{0\}$
and  first observe that
$\GO(\Phi^{R\theta})=\G(\Phi^{R\theta})$, as $\Phi^{R\theta}=0$ 
for $|x|<\rho$. Moreover, from Lemma \ref{phiprop} it results
$\G(\Phi^{R\theta})=\G(w)+o(\eps_{R})$, so that there exists
$R_{0}$ such that for every $R\geq R_{0}$, it results
$\Phi^{R\theta}\in A\cap A_{\Omega}$.
In order  to show \eqref{limtn},
let us consider $\theta\in \R^{N}\setminus\{0\},$ $R_{n}\to +\infty$ and set $T_{n}=T_{\Omega}(\Phi^{R_{n}\theta})$, which is well defined.
By  definition,  $T_{n}$ satisfies
\[
\|\Phi^{R_{n}\theta}\|_{\Omega}^{2}=\into 
\frac{f(T_{n}\Phi^{R_{n}\theta})}{T_{n}\Phi^{R_{n}\theta}}(\Phi^{R_{n}\theta})^{2}.
\]
Note that \eqref{defxi}, \eqref{limy}, \eqref{fmono}, \eqref{finfty} 
and \eqref{compact} imply
\begin{align*}
\|w^{R\theta}\|^{2}+o(\eps_{R})
=&
\intr \frac{f(T_{n}w^{R_{n}\theta})}{T_{n}w^{R_{n}\theta}}(w^{R_{n}\theta})^{2}
-\int_{B_{2\rho}(0)}\frac{f(T_{n}w^{R_{n}\theta})}{T_{n}w^{R_{n}\theta}}(w^{R_{n}\theta})^{2}
\\
&+\int_{\{\rho<|x|<2\rho\}}\frac{f(T_{n}\xi w^{R_{n}}\theta)}{T_{n}\xi w^{R_{n}}\theta}\left(\xi w^{R_{n}}\theta\right)^{2}
\\
=&
\intr \frac{f(T_{n}w^{R_{n}\theta})}{T_{n}w^{R_{n}\theta}}(w^{R_{n}\theta})^{2}
+o(\eps_{R}),
\end{align*}
Moreover, recalling \eqref{def:tras} and  performing a change of variable we obtain
\beq\label{equa:1}
\|w\|^{2}=
\intr \frac{f(T_{n}w)}{T_{n}w}w^{2}+o(\eps_{R}).
\eeq
Now, it is easy to see that $(T_{n})$ is bounded by contradiction,
since if it were not the case then one would obtain that
$\G(w)=0$ contradicting the fact that $w\in A$.
Then $(T_{n})$ has to be bounded, so that we can assume that,
up to a subsequence,  it converges to $T\in \R$, and $T\neq 0$ arguing 
again by contradiction.
Passing to the limit in \eqref{equa:1} we get
$$
\|w\|^{2}=\intr \frac{f(Tw)}{Tw}w^{2}.
$$
Then $T=1$ since $w\in {\mathcal N}$.
\edim

\bos\label{defpi}
Lemma \ref{tprop} provides the definition of the continuous projection
map
\[
\Pi_{\nhlo}: A_{\Omega}\to \nhlo,\quad \text{ as } \quad \Pi_{\nhlo}(u)=T_{\Omega}(u)u.
\]
\eos
Up to now, we have obtained asymptotic estimates on 
a function $\Phi^{R\theta}$.
In proving our existence results  we will use the following
convex combination of $\firy$ and $\firxz$, 
\beq \label{def:ur}
U_{t}^{R}=t\firy+(1-t)\firxz,
\quad  \text{with $x_{0}\in \R^{N},\,|x_{0}|=1$,
 $y\in \partial B_{2}(x_{0})$,  
 $t\in [0,1]$.}
\eeq
As a consequence, we will need  also
asymptotic informations on quantities 
involving  the functions $\Phi^{Ry}$ and $\Phi^{Rx_{0}}$
with $y\in \partial B_{2}(x_{0})$. 
Let us start with the following result.
\ble\label{scalprod}
Assume \eqref{fdergrowth}.
Let $x_{0}$ and $y$ be fixed in \eqref{def:ur}. Then it results
\[
\begin{split}
\left|\into \nabla \firy\nabla \firxz +\la\firy\firxz\right| 
\leq &
\frac12\intr f( w^{Ry}) w^{Rx_{0}}+\frac12\intr f( w^{Rx_{0}}) w^{Ry}
\\
& +o(\eps_{R}).
\end{split}
\]
\ele
\bdim
Notice that
\begin{align*}
\into \nabla \firy\nabla \firxz +\la\firy\firxz
=&
\into \nabla \xi\left[w^{Ry}\nabla \xi w^{Rx_{0}}+\xi(w^{Ry}\nabla w^{Rx_{0}}+\nabla w^{Ry}w^{Rx_{0}})
\right]
\\
&+\intr (\xi^{2}-1)\left[\nabla w^{Ry}\nabla w^{Rx_{0}}+
\la w^{Ry}w^{Rx_{0}}\right]
\\
&+\intr \nabla w^{Ry}\nabla w^{Rx_{0}}+
\la w^{Ry}w^{Rx_{0}}.
\end{align*}
As $w^{Ry}$ and $w^{Rx_{0}}$ are solutions of \eqref{Pinf}
we get
\[
\begin{split}
\into \nabla \firy\nabla \firxz +\la\firy\firxz
=&
\frac12\intr f( w^{Ry}) w^{Rx_{0}}+\frac12\intr f( w^{Rx_{0}}) w^{Ry}
\\
+\into \nabla \xi \left[w^{Ry}\nabla \xi w^{Rx_{0}}+
\right.&\left.\xi(w^{Ry}\nabla w^{Rx_{0}}+\nabla w^{Ry}w^{Rx_{0}})
\right]
\\
+\intr (\xi^{2}-1)\left[\nabla w^{Ry}
\nabla \right.&\left. 
w^{Rx_{0}}+\la w^{Ry}w^{Rx_{0}}\right].
\end{split}
\]
Then, the results is proved once one shows that 
the last two integrals on the right hand side are $o(\eps_{R})$.
First, we observe that \eqref{defxi} imply that there exists 
a positive constant $C_{0}$ such that
\begin{align*}
\left|\into \nabla \xi\left[w^{Ry}\nabla \xi w^{Rx_{0}}+\xi(w^{Ry}\nabla w^{Rx_{0}}+\nabla w^{Ry}w^{Rx_{0}})
\right] \right|
\leq &
C_{0}\int_{B_{2\rho}(0)}\big[(w^{Ry})^{2}+(w^{Rx_{0}})^{2}\big]
\\
+&C_{0}\int_{B_{2\rho}(0)}\big[|\nabla w^{Ry}|^{2}+|\nabla w^{Rx_{0}}|^{2}\big]
\\
\leq &
o(\eps_{R})
\end{align*}
where the last inequality is implied by \eqref{compact}.
In addition exploiting again \eqref{defxi} and \eqref{compact}
one obtains
\begin{align*}
\left|\intr (\xi^{2}-1)\left[\nabla w^{Ry}\nabla w^{Rx_{0}}+
\la w^{Ry}w^{Rx_{0}}\right]\right|
&\leq 
\frac12\int_{B_{2\rho}(0)}\left[|\nabla w^{Ry}|^{2}+
|\nabla w^{Rx_{0}}|^{2}\right]
\\
+\frac{\lambda}2 \int_{B_{2\rho}(0)}\left[(w^{Ry})^{2}+(w^{Rx_{0}})^{2}\right]
&\leq  o(\eps_{R})
\end{align*}
yielding  the conclusion.
\edim
In the following lemma we prove some crucial estimates
for the map $T_{\Omega}(U_{t}^{R})$.
\begin{lemma}\label{le:uniform}
Assume conditions \eqref{fprol}, \eqref{fdergrowth}, 
\eqref{fmono}, \eqref{finfty} and let $x_{0},\,y$ and 
$t$ be given in \eqref{def:ur}. Then,
there exists $R_{1}$ such that
 the following conclusions hold:
\\
i)\, $U^{R}_{t}\in A_{\Omega} $, for every $R\geq R_{1}$;
\\
ii)\, there exists  a positive constant $L$ such that
\beq\label{uniform}
|T_{\Omega}(U^{R}_{t})|\leq L, \quad
\text{for every $(t,R)\in [0,1]\times [R_{1},+\infty)$.}
\eeq
iii)\, For every $t_{0}\in (0,1)$, it  holds
\beq\label{lim:12}
\lim_{(t,R)\to (t_{0},+\infty)}T_{\Omega}(U^{R}_{t})=\dfrac1{t_{0}}\;\;
\text{if and only if}\;\; t_{0}=\frac12.
\eeq
\end{lemma}
\bdim
First of all, let us note that
\eqref{limy} and Lemmas \ref{lemma:eps}  and
\ref{scalprod} imply
that the following inequality holds for every 
$y\in \partial B_{2}(x_{0})$ (recall \eqref{def:ur})
\beq\label{normat}
\begin{split}
\|U^{R}_{t}\|_{\Omega}^{2}
=&
\|U^{R}_{t}\|^{2}
=
t^{2}\|\firy\|^{2}+(1-t)^{2}\|\firxz\|^{2}
\\
+& 2t(1-t)\intr \nabla\firy\nabla \firxz +\firy\firxz
\\
\leq& t^{2}\|w^{Ry}\|+(1-t)^{2}\|w^{Rx_{0}}\|^{2}+2t(1-t)\eps_{R}+o(\eps_{R}).
\end{split}
\eeq
Then, recalling the definition of $\GO$ stated in Lemma \ref{tprop},
using again \eqref{limy}, and applying Lemma \ref{scalprod},  one 
obtains
\begin{align*}
\GO(U^{R}_{t})\leq 
& t^{2}\|w^{Ry}\|^{2}+(1-t)^{2}\|w^{Rx_{0}}\|^{2}
-l_{\infty}t^{2}\|w^{Ry}\|^{2}-l_{\infty}(1-t)^{2}\|w^{Rx_{0}}\|^{2}
\\
&+2t(1-t)\eps_{R}+2l_{\infty}t(1-t)\intr w^{Ry}w^{Rx_{0}}+o(\eps_{R})
\\
=& t^{2}\G (w)+(1-t)^{2}\G(w)+2t(1-t)\eps_{R}
\\
&+2l_{\infty}t(1-t)\intr w^{Ry}w^{Rx_{0}}+o(\eps_{R}).
\end{align*}
Thanks to \eqref{behavior:w} we can apply Lemma 2.1 in
\cite{acclpa} with $\overline{\mu}=\sqrt{\la}/2$ to deduce
\[
\lim_{R\to+\infty}\intr w^{Ry}w^{Rx_{0}}\leq
\lim_{R\to+\infty}e^{-\sqrt{\la}R}=0.
\]
Then, as $\G(w)<0$,  there exists $R_{1}$ such that for every 
$R>R_{1}$ it holds \(\GO(U^{R}_{t})<0\), for every 
$y\in \partial B_{2}(x_{0})$, yielding conclusion {\it i)}.
In order to prove conclusion {\it ii)}, first note that
$T_{\Omega}(U^{R}_{t})$ is well defined; then,
we argue by contradiction and suppose that there exist
sequences $R_{n}\to +\infty $, 
$t_{n}\in [0,1]$ and $y_{n}\in \partial B_{2}(x_{0})$ 
such that
\beq\label{tdiv}
U_{n}:=t_{n}\phiyn+ (1-t_{n})\firxzn\;\text{satisfies}\;
T_{n}:=T_{\Omega}(U_{n})\to +\infty.
\eeq
Since $t_{n}\in [0,1]$ we can suppose, up to a subsequence, that
there exists $t_{0}\in[0,1]$ such that $t_{n}\to t_{0}$.
Three cases may occur:
either  $t_{0}=0$, or $t_{0}=1$ or $t_{0}\in (0,1)$.
\\
Suppose first that $t_{0}=0$, and observe that
$T_{n}U_{n}\in {\mathcal N}_{\Omega}$
so that
\beq\label{eq:t}
\|U_{n}\|_{\Omega}^{2}=
\into
\frac{f(T_{n}U_{n})}{T_{n}U_{n}}U_{n}^{2}.
\eeq
Then, taking into account  that $t_{n}\to 0,\,R_{n}\to+\infty$
and exploiting \eqref{normat}, 
we obtain
\begin{align}\label{lim1}
\lim_{n\to+\infty}
\|U_{n}\|_{\Omega}^{2}
=\|w\|^{2}.
\end{align}
With respect to the right hand side of \eqref{eq:t} we observe that
the property of the function $\xi$, joint with
assumptions  \eqref{fmono} and \eqref{finfty}   give
that there exists a positive constant $C_{1}$ such that
\beq\label{limf0}
\begin{split}
\into
\frac{f(T_{n}U_{n})}{T_{n}U_{n}}
\left[t^{2}_{n}(\phiyn)^{2}
+2t_{n}(1-t_{n}) \phiyn\firxzn\right]
&\leq 
C_{1}t_{n}\|w\|^{2}_{2}.
\end{split}
\eeq
Moreover, one observes that
\[
\into
\frac{f(T_{n}U_{n})}{T_{n}U_{n}}(\firxzn)^{2}
=
\into
\frac{f(T_{n}U_{n})}{T_{n}U_{n}}
(w^{R_{n}x_{0}})^{2}
+
\into\frac{f(T_{n}U_{n})}{T_{n}U_{n}}\left[(\firxzn)^{2}-(w^{R_{n}x_{0}})^{2}\right]
\]
and assumptions \eqref{fmono}, \eqref{finfty} and \eqref{limy}
yield
\beq\label{limf1}
\into
\frac{f(T_{n}U_{n})}{T_{n}U_{n}}\left[(\firxzn)^{2}-(w^{R_{n}x_{0}})^{2}\right]
\leq
2l_{\infty}\int_{B_{2\rho}(0)} (w^{R_{n}x_{0}})^{2}\leq 
o(\eps_{R_{n}}).
\eeq
In addition, since
\beq\label{f:tras}
\begin{split}
\into 
\frac{f(T_{n}U_{n})}{T_{n}U_{n}}(w^{R_{n}x_{0}})^{2}
=
\intr 
\frac{f(T_{n}(t_{n}\xi w^{R_{n}y_{n}}+(1-t_{n})\xi w^{R_{n}x_{0}}))}
{T_{n}(t_{n}\xi w^{R_{n}y_{n}}+(1-t_{n})\xi w^{R_{n}x_{0}})}(w^{R_{n}x_{0}})^{2}
\\
=
\intr
\frac{f(T_{n}(t_{n}\xi^{-R_{n}x_{0}} w^{-R_{n}(y_{n}+x_{0})}+(1-t_{n})
\xi^{-R_{n}x_{0}} w))}
{T_{n}(t_{n}\xi^{-R_{n}x_{0}} w^{-R_{n}(y_{n}+x_{0})}+(1-t_{n})
\xi^{-R_{n}x_{0}} w)}w^{2}
\end{split}
\eeq
and as  \eqref{defxi} and \eqref{tdiv} imply
$$
T_{n}(t_{n}\xi^{-R_{n}x_{0}} w^{-R_{n}(y_{n}+x_{0})}+(1-t_{n})
\xi^{-R_{n}x_{0}} w)
\geq
T_{n}(1-t_{n})\xi(x+R_{n}x_{0}) w(x)\to +\infty
$$
almost  everywhere in $\R^{N}$,
\eqref{finfty} and \eqref{fmono} yields
\beq\label{limf2}
\lim_{n\to+\infty}\into 
\frac{f(T_{n}U_{n})}{T_{n}U_{n}}(w^{R_{n}x_{0}})^{2}
=l_{\infty}\intr w^{2}.
\eeq
Then \eqref{lim1}, \eqref{limf0}, \eqref{limf1}, \eqref{limf2}
allow us to pass to the list in \eqref{eq:t} to obtain
$\G(w)=0$,  which is a contradiction as $w\in A$ (defined in Lemma \ref{tprop}).
\\
The case in which $t_{0}=1$  is similar to the case $t_{0}=0$ 
by exchanging the role of  $t_{n}$ with the one of $1-t_{n}$.
Let us handle the third case and  suppose that $t_{0}\in (0,1)$ and come back to \eqref{normat} 
to observe that \eqref{limy} and
Lemma  \ref{scalprod} 
imply that
\beq\label{unmezzo}
\lim_{n\to+\infty}\|U_{n}\|_{\Omega}^{2}=
\left[t^{2}_{0}+(1-t_{0})^{2}\right]
\|w\|^{2}.
\eeq
In order to study the right hand side of \eqref{eq:t} notice that, by
\eqref{fmono} and \eqref{defxi}
\[
\begin{split}
0\leq \limsup_{n\to+\infty}
\into\frac{f(T_{n}U_{n})}{T_{n}U_{n}}\phiyn\firxzn
\leq l_{\infty}\lim_{n\to+\infty} \intr w(x)w(x-R_{n}(x_{0}-y))
=0
\end{split}
\]
and as $t_{0}\in (0,1)$ it results
\[
\lim_{n\to +\infty}
2t_{n}(1-t_{n})\into\frac{f(T_{n}U_{n})}{T_{n}U_{n}}\phiyn\firxzn
=0.
\]
Then, taking into account that
\eqref{limf1} holds also for $y$ in the place of $x_{0}$,
one obtains
\[
\begin{split}
\lim_{n\to+\infty}\into\frac{f(T_{n}U_{n})}{T_{n}U_{n}}
U_{n}^{2}
=&
\lim_{n\to+\infty}t^{2}_{n}
\into\frac{f(T_{n}U_{n})}{T_{n}U_{n}}
(w^{R_{n}y_{n}})^{2}
\\
&+
\lim_{n\to+\infty}(1-t_{n})^{2}\into\frac{f(T_{n}U_{n})}{T_{n}U_{n}}(w^{R_{n}x_{0}})^{2}.
\end{split}
\]
Arguing as in \eqref{f:tras} in the above integrals separately,
exploiting  \eqref{defxi} and \eqref{tdiv} and
taking into account that $ t_{0}\in(0,1)$
  one deduces
$$
T_{n}\left[t_{n}\xi^{-R_{n}y_{n}}w+(1-t_{n})\xi^{-R_{n}y_{n}}w^{R_{n}(x_{0}-y_{n})} \right]\geq T_{n}t_{n}\xi(x-R_{n}y_{n})w\to+\infty \quad \text{a.e. in} \,\Omega.
$$
Therefore, we can use \eqref{finfty} and \eqref{fmono}  to apply Lebesgue dominated convergence Theorem and obtain
\beq\label{unmezzo2}
\lim_{n\to+\infty}\into\frac{f(T_{n}U_{n})}{T_{n}U_{n}}
U_{n}^{2}
=\left[t_{0}^{2}+(1-t_{0})^{2}\right]l_{\infty}\|w\|_{2}^{2}.
\eeq
Finally, passing to the limit in \eqref{eq:t} and using 
 \eqref{unmezzo}, \eqref{unmezzo2} it follows
\beq\label{unmezzo3}
\left[t_{0}^{2}+(1-t_{0})^{2}\right]\|w\|^{2}
=\left[t_{0}^{2}+(1-t_{0})^{2}\right]l_{\infty}\|w\|_{2}^{2},
\eeq
which again contradicts the fact that $w\in A$, yielding 
\eqref{uniform}.
In order to prove \eqref{lim:12}, let us first show that 
\beq\label{limmezzo}
\lim_{(t,R)\to (\frac12,+\infty)}T_{\Omega}(U^{R}_{t})=2.
\eeq
Arguing again by contradiction and
supposing that there exist $\delta>0$ and  
subsequences $R_{k}\to+\infty$,
$ t_{k} \to 1/2$ and $y_{k}\in \partial B_{2}(x_{0})$ 
such that, the sequence 
$$
T_{k}:=T_{\Omega}\left(t_{k}\Phi^{R_{k}y_{k}}+(1-t_{k})\Phi^{R_{k}x_{0}}\right)
\; \;\text{satisfies} \;\; |T_{k}-2|\geq \delta.
$$
Estimate \eqref{uniform} implies that there exists 
$T$ such that, up to a subsequence, $T_{k}\to T$.
Then \eqref{eq:t} becomes
\beq\label{eq:t2}
\|U_{k}\|_{\Omega}^{2}=
\into
\frac{f(T_{k}U_{k})}{T_{k}U_{k}}U_{k}^{2}
\eeq
and \eqref{unmezzo} allows to take the limit on the left hand side;
on the other hand, we can argue in analogous way as in the case 
$t_{0}\in (0,1)$  to obtain
\[
\frac12\|w\|^{2}=
\intr \frac{f\left(\frac{T}2 w\right)}{\frac{T}2 w}\frac{w^{2}}2.
\]
As $w$ is a solution of the limit problem \eqref{Pinf} this implies
\[
\intr \left[\frac{f(w)}{w}-\frac{f\left(\frac{T}2 w\right)}{\frac{T}2 w}
\right]w^{2}=0,
\]
yielding \eqref{limmezzo} thanks to \eqref{fmono}.
On the other hand, suppose that  
there exists $t_{0}\in (0,1)$ such that the 
limit in \eqref{lim:12} holds and
use  \eqref{limy} and Lemma \eqref{scalprod} to take the limit
in \eqref{eq:t2} and obtain
$$
t_{0}^{2}\|w\|^{2}+(1-t_{0})^{2}\|w\|^{2}=
t_{0}^{2}\intr \frac{f(w)}{w}w^{2}+
\intr \frac{f(w(1-t_{0})/t_{0})}{w(1-t_{0})/t_{0}}
(1-t_{0})^{2}w^{2}
$$
that is 
\[
0
=
(1-t_{0})^{2} \intr \left[\frac{f(w)}w-
\frac{f(w(1-t_{0})/t_{0})}{w(1-t_{0})/t_{0}}
\right]w^{2}.
\]
Finally, \eqref{fmono} implies that $t_{0}=1/2$ yielding the conclusion.
\edim

\ble\label{ack}
Assume condition \eqref{fdergrowth} and let $x_{0}$,
$y$ and $t$ be given in \eqref{def:ur}.
Then, it results
\[
\begin{split}
\left|
\intr 
\big[
F(T_{R}^{t}U^{R}_{t})
\right.
-F(T_{R}^{t}t\firy)&-F(T_{R}^{t}(1-t)\firxz)
\big]
\\
+\intr 
\big[
-f(T_{R}^{t}t\firy)T_{R}^{t}(1-t)\firxz
&
-f(T_{R}^{t}(1-t)\firxz)T^{R}_{t}t\firy
\big]\Big|
\leq o(\eps_{R}),
\end{split}
\]
where  $T^{R}_{t}=T_{\Omega}(U^{R}_{t})$.
\ele
\bdim
We use Lemma 2.2 in \cite{acclpa}, with $\alpha:= \min \{(p_1 +1)/4, 1\}$ and $p_1$ as in (\ref{fdergrowth}), and
we take into account \eqref{defxi} and  \eqref{uniform} 
to obtain that
\[\begin{split}
\left|\intr 
\big[
F(T_{R}^{t}U^{R}_{t})
-F(T_{R}^{t}t\firy)-F(T_{R}^{t}(1-t)\firxz)
\big]\right.&
\\
+\intr 
\big[
-f(T_{R}^{t}t\firy)T_{R}^{t}(1-t)\firxz
-f(T_{R}^{t}(1-t)\firxz)T_{R}^{t}t\firy
\big] 
\Big| &
\\
\leq 
\left[T_{R}^{t}t\right]^{2\alpha}\left[T_{R}^{t}(1-t)\right]^{2\alpha} 
\intr |\firy|^{2\alpha}
|\firxz|^{2\alpha}&
\\
\leq L^{4\alpha }\intr   w^{2\alpha}(x) w^{2\alpha}(x-R(y-x_{0})) \,dx.\quad\,&
\end{split}\]
Since $2\alpha>1$, we have that, denoting with 
$\vfi(x)=w^{2\alpha}(x)$,  the first  hypothesis 
in Lemma \ref{bali} is satisfied with $c=0$.
Then Lemma \ref{bali} implies that
$$
\lim_{R\to \infty }(2R)^{(N-1)/2}e^{2R\sqrt{\lambda}}
\intr  w^{2\alpha}(x)w^{2\alpha}(x-R(y-x_{0}))\,dx=0.
$$
Then, Lemma \ref{lemma:eps} yields the conclusion.
\edim
\ble\label{le:lower}
Assume \eqref{fprol}, \eqref{fdergrowth}, \eqref{fmono}
and let $x_{0}$ and $y$ be given in \eqref{def:ur}. Then,
There exist $C=C(\|w\|_{\infty})>0$
such that, for all  $\tau_1, \tau_2\in [0,+\infty) $, 
it results
\[
\dys\intr f(\tau_1 w^{Rx_{0}})\tau_2 w^{Ry}\geq \min\{\tau_{1},\tau_{2}\} O(\eps_{R}),\,
\dys\intr f(\tau_1 w^{Ry})\tau_2 w^{Rx_{0}}\geq  \min\{\tau_{1},\tau_{2}\} O(\eps_{R}).
\]
\ele
\bdim
As, performing a change of variable we have
\[
\begin{split}
\intr f(\tau_1 w^{Rx_{0}})\tau_2 w^{Ry}=
\intr f(\tau_1 w)\tau_2 w^{R(x_{0}-y)}, 
\\
\intr f(\tau_1 w^{Ry})\tau_2 w^{Rx_{0}}=
\intr f(\tau_1 w)\tau_2 w^{R(y-x_{0})}, 
\end{split}
\]
it is enough to show the first inequality,
since the other will follow by a similar argument.
In order to do this, taking into account Lemma \ref{lemma:eps},
it is sufficient to show that there exists a constant $C$ such that
$$
\intr f(\tau_1 w^{Rx_{0}})\tau_2 w^{Ry}\geq C R^{-(N-1)/2}
e^{-2R\sqrt{\la}}.
$$
Note that, the positive minimum of $w(x)$ 
in the ball $B_{1}(0)$ is achieved for $|x|=1$; let us denote
this minimum value as $\alpha_{0}>0$.
As a consequence of hypothesis \eqref{fmono} we
have that the function $g(u)=f(u)/u$ is continuous and monotone 
increasing in the interval $[\alpha_{0},\|w\|_{\infty}]$, so that
 $g(u)\geq C :=g(\alpha_{0})$. Then we obtain
\begin{align}\label{fstima}
\intr f(\tau_1w^{Rx_{0}})\tau_2w^{Ry}&
\geq \tau_1 \tau_2\int_{B_{1}(0)} \frac{f(\tau_1w)}{\tau_1w}ww^{R(x_{0}-y)} 
\\
\nonumber
&\geq  \tau_{1}\tau_{2} C \int_{B_{1}(0)}w(x)w(x-R(x_{0}-y))\,dx.
\end{align}
Note that, for every $x\in B_{1}(0)$ it results for every $R\geq1$
\beq\label{stima}
 2R - 1\leq |R(y-x_0)| - |x| \leq | x-R(y-x_0) |\leq 2R+1\leq 3R.
\eeq
In the following we will denote with $C$ possibly different 
positive constants.
Condition \eqref{behavior:w} implies that for sufficiently large $R$
there exists a constant $C$ such that
$$
w(x-R(x_{0}-y))\geq C |x - R(y-x_0)|^{-(N-1)/2}e^{-\sqrt{\lambda}|x - R(y-x_0)|}
$$
and \eqref{stima} gives 
$$
w(x-R(x_{0}-y))\geq C R^{-(N-1)/2}e^{-2R\sqrt{\la}}.
$$
Using this inequality in \eqref{fstima} and applying 
Lemma \ref{lemma:eps} yield the conclusion.
\edim
\begin{lemma} \label{tvm}
Assume \eqref{fprol}, \eqref{fdergrowth}, \eqref{fmono},  
and let $x_{0}$ and $y$ be given in \eqref{def:ur}.
Then, for every $\tau\in [0,+\infty)$, 
it results
\begin{align}  
\label{stima:f}
\left|\intr f(\tau \Phi^{Ry})\Phi^{Rx_{0}}-f(\tau w^{Ry})w^{Rx_{0}}\right|
&\leq 2D \max\{\tau ^{p_{1}},\tau ^{p_{2}}\} o(\eps_{R}) ,
\\
\label{inttvm}
\left|\int_{\mathbb{R}^N} \tau f(w^{Rx_{0}}) w^{Ry}
- f(\tau w^{Rx_{0}}) w^{Ry}\; \right|
&\leq |\tau -1| \left[1+D\tau^{\max\{p_{1},p_{2}\}-1}\right] O (\varepsilon_R).
\end{align}
\end{lemma}
\bdim
The first inequality is a direct consequence of \eqref{fmono}, \eqref{fdergrowth},
\eqref{defxi} and \eqref{compact}. Indeed it results
\begin{align*}
\intr \left|f(\tau\firy)\firxz -f(\tau w^{Ry})w^{Rx_{0}}
\right| 
\leq 
2\int_{B_{2\rho}(0)} & f(\tau  w^{Ry})w^{Rx_{0}}
\\
\leq 
2D \max\{\tau ^{p_{1}},\tau ^{p_{2}}\}\int_{B_{2\rho}(0)}
[(w^{Ry})^{p_{1}}+(w^{Ry})^{p_{2}}]w^{Rx_{0}}
\leq 
2D &\max\{\tau ^{p_{1}},\tau ^{p_{2}}\}o(\eps_{R}).
\end{align*}
In order to show \eqref{inttvm}, consider the function $
g(\tau):=\tau f(w^{Rx_{0}}) - f(\tau w^{Rx_{0}}),$
and suppose,  without loss of generality, that $ \tau>1$. 
The mean value theorem implies that there exists 
$\theta\in (1,\tau)$ such that
\begin{align*}
g(\tau)=g(\tau)-g(1)&= g'(\theta ) (\tau -1)= \left[f(w^{Rx_{0}}) - f'(\theta w^{Rx_{0}}) w^{Rx_{0}}\right] (\tau -1).
\end{align*}
Substituting in the integral in (\ref{inttvm}) we obtain
\begin{align*}
\left|
\int_{\mathbb{R}^N} [\tau f(w^{Rx_0})- f(\tau w^{Rx_0})] 
w^{Ry}\;  \right|
=&
|\tau -1|  {\Big |}\int_{\mathbb{R}^N} 
\left[f(w^{Rx_0})- f'(\theta w^{Rx_0})w^{Rx_{0}}\right] w^{Ry}\; 
{\Big |}
\\
\leq &
|\tau -1|\int_{\mathbb{R}^N} f(w^{Rx_0}) w^{Ry}+\left| f'(\theta w^{Rx_0})\right|w^{Rx_0} w^{Ry}.
\end{align*}
Applying  Lemma \ref{lemma:eps}, hypothesis
\eqref{fdergrowth} ad taking into account that $\theta\in(1,\tau)$
one obtains
\begin{align*}
\left|
\intr [\tau f(w^{Rx_0})- f(\tau w^{Rx_0})] w^{Ry}\;  \right|
&\leq 
|\tau-1|\eps_{R} 
\\
+
D|\tau-1| \tau^{\max\{p_{1},p_{2}\}-1}
&\intr \left[
( w^{Rx_0})^{p_1}w^{Ry}\; + (w^{Rx_0})^{p_2}w^{Ry}\; 
\right]  .
\end{align*}
Applying Lemma \ref{bali}
with $\vfi_{1}=w$ and $\vfi_{2}=w^{p_{1}}$ 
(see also the argument of the proof of  Lemma \ref{lemma:eps}),
we obtain that 
$$
\int_{\mathbb{R}^N}  
\left[ (w^{Rx_0})^{p_1}w^{Ry}\;+ 
(w^{Rx_0})^{p_2}w^{Ry}\; \right]  \leq
O(\varepsilon_R)\;
$$
so that  the conclusion follows.
\edim
\section{Compactness Condition}\label{compactness}
In this section we will find the level interval where Palais-Smale condition
holds.
\begin{lemma}\label{equal}
Assume \eqref{fprol}, \eqref{fmono},\eqref{finfty}. There holds
\[
m_{\Omega}=m
\]
and $m_{ \Omega}$ is not attained.
\end{lemma}
\bdim
Since any $u \in H^1_0(\Omega)$ can be extended as zero outside $\Omega$, we may consider $H^1_0(\Omega) \subset H^1(\mathbb{R}^N)$ and so $m_{\Omega} \geq m\;.$
On the other hand, 
applying Lemma \ref{tprop} to the sequence 
$\phi_n:=\Phi^{R_{n}\theta}$, with $R_{n}\to +\infty$,
it follows that  there exists $n_{0}$ such that for every $n\geq n_{0}$
there exists $T_n >0$ 
such that $T_n (\phi_{n})\phi_{n} \in \mathcal{N}_{\Omega}.$
Then, taking into account \eqref{m.omega},  \eqref{limyI} and \eqref{limtn} yield
\[
m_{\Omega} \leq I_{\Omega}(T_n (\phi_{n})\phi_{n})=
I_{\Omega}(T_n (\Phi^{R_{n}\theta})\Phi^{R_{n}\theta})\to 
I(w)=m\,.
\]
Finally, assume that there exists $u\in H^{1}_{0}(\Omega)$
such that $I_{\Omega}(u)=m_{\Omega}$. Then
$\ov{u}$, extension by zero of $u$ outside $\Omega$, would 
be a minimizer of the problem at infinity not positive on
the whole $\R^{N}$ contradicting the maximum principle.
\edim
The remaining of this section is devoted to recover compactness
properties for $I_{\Omega}$. The main difficulty in our 
context is due to the fact that $f$ is asymptotically linear.
Moreover, we will need a compactness property just for 
sequences belonging to the Nehari manifold $\nhlo$, 
then we  first  prove the following result.

\begin{lemma}\label{ceramibounded}
Assume \eqref{fprol}, \eqref{fdergrowth}, \eqref{fmono},
\eqref{finfty}, \eqref{NQ}, \eqref{iponh}.
Let $d\in \R^{+}$  and $(u_n)$ be such that
\beq\label{cesequence}
(u_n)\in\nhlo \quad\text{and}\quad \lim_{n\to+\infty} I_\Omega(u_n)= d>0  
\eeq
then $(u_{n})$ is bounded. 
\end{lemma}
\bdim
The proof can be started as in Proposition 3.20 in \cite{mamope} 
and it can be concluded as in Lemma 5.3 in \cite{mamiso}. 
We give here a brief summary.
First, note that as $u_{n}\in \nhlo$ it results
\beq\label{equaun}
I_{\Omega}(pu_{n})=\frac{p^{2}}{2}\|u_{n}\|_{\Omega}^{2}
-\into F(pu_{n})dx
=\into \frac{p^{2}}{2}f(u_{n})u_{n}-F(pu_{n}).
\eeq
For every $v\in \nhlo $, let us consider the one-variable function $h:\R^{+}\to \R$ 
\[
h_{v}(p):= \frac{p^{2}}{2}f(v)v-F(pv),\quad I(pv)=\into h_{v}(p).
\]
Notice that \eqref{cesequence} implies that there exists a positive 
constant $C=C(d)$ such that
\[
|\io(u_{n})|\leq C,
\]
so that, condition \eqref{fmono} implies that, for every
fixed $x\in \R^{N}$,  $h_{u_{n}}(p)$ has a unique  maximum point 
for $p=1$ so that
\(\io(pu_{n})\leq \io(u_{n})\leq C\).
Recalling \eqref{equaun}, one obtains
$$
\int_{\Omega} F(pu_{n})\geq \frac{p^{2}}2\|u_{n}\|^{2}-C.
$$
Then, assuming by contradiction that,
up to a subsequence,  $\|u_{n}\|_{\Omega}\to +\infty$
and setting $v_{n}=p_{n}u_{n}$ with $p_{n}=2\sqrt{C}/\|u_{n}\|_{\Omega}$,
we obtain the uniform lower bound
\[
\into F(v_{n})\geq C.
\]
This lower bound, hypothesis \eqref{fdergrowth}
and Lions Lemma (see \cite{lions}) imply that there exist positive numbers $r$ and $\delta$ and a sequence $(y_{n})\in \R^{N}$ 
such that
\beq\label{boundlions}
\liminf_{n\to\infty}\int_{B_{r}(y_{n})}v^{2}_{n}dx\geq \delta.
\eeq
Then, we have to handle two different possible situations:
either  $y_{n}$ is bounded, or, up to a subsequence, 
$|y_{n}|\to+\infty$ as $n\to +\infty$.
\\
In the first case we deduce from \eqref{boundlions} and recalling that $v_{n}\equiv 0$ outside $\Omega$, that there exists
$r_{1}>\rho$ such that 
$$
\liminf_{n\to\infty}\int_{B_{r_{1}}(0)}v^{2}_{n}dx\geq \frac{\delta}2
$$
and the same lower bound holds for the weak limit $v$ of $v_{n}$
(up to a subsequence). Then, as $v_{n}\equiv 0$ 
in $\R^{N}\setminus\Omega$ there exists a subset 
$\Lambda\in B_{r_{1}}(0)\cap \Omega$ 
with positive measure and 
such that $v(x)>0$ in $\Lambda$, so as 
$u_{n}(x)=2\sqrt{C}v_{n}(x)\|u_{n}\|_{\Omega}$,
it follows that
$u_{n}(x)\to+\infty$ for all $x\in \Lambda$.
This immediately leads to a contradiction using \eqref{fprol} 
\eqref{NQ} and  \eqref{cesequence} (for more details see Lemma 5.3 in
\cite{mamiso}).
\\
Then $(y_{n})$ cannot be bounded, and, up to a subsequence, 
we obtain $|y_{n}|\to+\infty$; moreover,
we can assume that $B_{r}(y_{n})\subset \Omega$
for $n$ sufficiently large. Then, it is possible to argue 
analogously on the sequence $\tilde{v}_{n}=v_{n}(x+y_{n})$, 
obtaining, as before, a set  $\Lambda\subseteq B_{r}(0)$ 
with positive measure such that $u_{n}(x+y_{n})\to+\infty$ 
for every $x\in \Lambda$. 
The contradiction follows again by using
\eqref{NQ} (see Lemma 5.3 in \cite{mamiso}). Indeed, 
\eqref{cesequence} yields
\begin{align*}
C &\geq\limsup_{n\to+\infty}\io(u_{n})-\langle \io'(u_{n}),u_{n}\rangle
\geq \int_{B_{r}(y_{n})}\frac12 f(u_{n})u_{n}-F(u_{n})
\\
&=\int_{B_{r}(0)}
\frac12 f(u_{n}(x+y_{n}))u_{n}(x+y_{n})-F(u_{n}(x+y_{n}))dx
\\
&\geq \int_{\Lambda}
\frac12 f(u_{n}(x+y_{n}))u_{n}(x+y_{n})-F(u_{n}(x+y_{n}))dx
\end{align*}
and the last integral goes to plus infinity thanks to \eqref{NQ}.
Then, we reach a contradiction, proving the Lemma.
\edim
\bos
Notice that in order to prove that $(u_{n})$ is bounded
the classical information $I_{\Omega}'(u_{n})\to 0$
can be substituted by the  information
$(u_{n})\subseteq \nhlo$.
Moreover, note that, in order to prove Lemma \ref{ceramibounded},
it is sufficient to assume \eqref{fdergrowth} for $k=0$.
\eos
In the following lemma we will show that $(u_{n})$ is a Palais-Smale
sequence in the whole space. 
\ble\label{nuovo}
Assume \eqref{fprol}, \eqref{fdergrowth}, \eqref{fmono}, 
\eqref{finfty}, \eqref{NQ}, \eqref{iponh}.
Let $(u_{n})$ satisfy  \eqref{cesequence} and be such that
\beq\label{ipo} 
\nabla_{\nhlo} I_{\Omega}(u_n) \to 0 .
\eeq
Then $\io'(u_{n})\to 0$ in $H^{-1}(\Omega)$.
\ele
\bdim
First of all, 
from \eqref{ipo}, recalling \eqref{nehari}, 
we obtain a sequence $(\mu_{n})\subseteq\R$ 
such that 
\beq\label{equavinc}
\io'(u_{n})-\mu_{n}N_{\Omega}'(u_{n})\to 0,
\eeq
where $N_{\Omega}$ is defined in Remark \ref{natu}.
Moreover,   we can use Lemma \ref{ceramibounded}
to obtain that $(u_{n})$ is bounded, so that
there exists $u\in H^{1}_{0}(\Omega)$
such that, up to a subsequence, $u_{n}\rightharpoonup u$ 
weakly in $\homega$ and $u_{n}\to u$ almost everywhere.
Let us first show that $(\mu_{n})$ is bounded arguing by contradiction,
so that we assume that, up to a subsequence, $|\mu_{n}|\to+\infty$ and
set $t_{n}=\mu_{n}/|\mu_{n}|$. Since $|t_{n}|=1$ there
exists $t_{0}$ such that, up to a subsequence, $|t_{0}|=1$ and $t_{n}\to t_{0}$ in $\R$.
Moreover,  it results
$$
\frac1{|\mu_{n}|}\io'(u_{n})-t_{n} N_{\Omega}'(u_{n})\to 0.
$$
Since $|\mu_{n}|\to +\infty$ and as
$(u_{n})$ is  bounded \eqref{fdergrowth} implies
\[
N'_{\Omega}(u_{n})\to 0,
\]
because  $t_{n}\to t_{0}\neq 0$.
Then, using \eqref{fprol}, and taking into account that 
$u_{n}\in \nhlo$  we have
\beq\label{nprime}
\begin{split}
0 &=\lim_{n\to+\infty} \langle N_{\Omega}'(u_{n}),u_{n}\rangle=
\lim_{n\to+\infty}\into \left[f(u_{n})u_{n}-f'(u_{n})u_{n}^{2}\right]
\\
&=\lim_{n\to+\infty}\into \left[f(u_{n}^{+})u_{n}^{+}-f'(u_{n}^{+})(u_{n}^{+})^{2}\right].
\end{split}
\eeq
Assume that there exists $\delta>0$ such that
$$
\sup_{y\in \R^{N}}\int_{B_{r}(y)}|u^{+}_{n}|^{2} \geq \delta
$$ 
then there exists $(y_{n})$ such that 
\beq\label{intoy}
\int_{B_{r}(y_{n})}|u^{+}_{n}|^{2} \geq \delta
\eeq
arguing as in Lemma \ref{ceramibounded}, we observe that 
two  cases may occur, either $|y_{n}|$ is bounded,
or, up to a subsequence, $|y_{n}|\to +\infty$.
In the first case, there exists $r_{1}>0$ such that 
$\|u^{+}_{n}\|^{2}_{L^{2}(B_{r_{1}}(0))}\geq \delta$.
Then, since $u^{+}_{n}$ strongly converges to some $v$,
in  $L^{2}_{{\rm loc}}(\R^{N})$, it follows that
$v\geq 0$, and $v$ satisfies
$\|v\|^{2}_{L^{2}(B_{r_{1}}(0))}\geq \delta$.
Moreover, recalling that $u^{+}_{n}\equiv 0$ 
outside $\Omega$, we deduce that there exists
a measurable set $\Lambda\subseteq \Omega\cap B_{r_{1}}(0)$ with 
positive measure, such that $v(x)>0$ for $x\in \Lambda$.  Then, using \eqref{NQ}, \eqref{iponh} and applying
Fatou Lemma, one obtains
\beq\label{fine1}
\begin{split}
0&=\lim_{n\to+\infty}-
\langle N'_{\Omega} (u_{n}),u_{n}\rangle 
=
\lim_{n\to+\infty}
\into f'(u^{+}_{n})(u^{+}_{n})^{2}-f(u^{+}_{n})u^{+}_{n}
\\
&\geq \lim_{n\to+\infty} \int_{\Lambda}
f'(u^{+}_{n})(u^{+}_{n})^{2}-f(u^{+}_{n})u^{+}_{n}
=\int_{\Lambda}f'(v)v^{2}-f(v)v>0,
\end{split}
\eeq
that is a contradiction, so that
 $|y_{n}| $ cannot be bounded, and
 $|y_{n}|\to+\infty$, up to a subsequence.
Then, we define $\hu_{n}(\cdot):=u^{+}_{n}(\cdot+y_{n})$ and observe that
there exists $\hu\geq 0$ almost everywhere in $\R^{N}$,   
such that $\hu_{n}\rightharpoonup \hu$ 
weakly in  $H^{1}_{0}(\Omega)$, 
strongly in $L^{p}_{{\rm loc}}(\R^{N})$ 
and almost everywhere. 
From the strong convergence in $L^{2}_{\rm loc}(\R^{N})$ and
using \eqref{intoy},
one deduces that $\|\hu\|^{2}_{L^{2}(B_{r_{1}}(0))}\geq \delta$.
Then there exists a measurable set $\Lambda\subseteq B_{r}(0)$ 
with positive measure, and such that $\hu(x)> 0$ for every $x\in \Lambda$.
Moreover, as $|y_{n}|\to \infty $, we can assume that 
$B_{r}(y_{n})\subseteq \Omega$.
Then, arguing as in  \eqref{fine1} we get
\begin{align*}
0&=\lim_{n\to+\infty}-
\langle N'_{\Omega} (u_{n}),u_{n}\rangle 
\geq \lim_{n\to+\infty} \int_{B_{r}(y_{n})}
f'(u^{+}_{n})(u^{+}_{n})^{2}-f(u^{+}_{n})u^{+}_{n}
\\
&=\lim_{n\to+\infty}\int_{B_{r}(0)}f'(\hu_{n})\hu_{n}^{2}-f(\hu_{n})\hu_{n}
\geq \int_{\Lambda}f'(\hu)\hu-f(\hu)\hu>0,
\end{align*}
so that \eqref{iponh} yields again a contradiction, showing that
\[
\lim_{n\to+\infty}\sup_{y\in \R^{N}}
\int_{B_{r}(y)}|u^{+}_{n}|^{2}=0.
\]
Then, Lions' Lemma  \cite{lions} implies
$$
u^{+}_n\rightarrow 0 \quad \mbox{in} \ L^{p}(\mathbb{R}^N), \ \mbox{for} \ \mbox{any} \ 2<p<2^*
$$
and \eqref{fdergrowth} gives
\[
\lim_{n \to \infty} \int_\Omega f(u^{+}_{n})u^{+}_{n} \to 0\;, \quad 
\lim_{n \to \infty} \int_\Omega F(u^{+}_{n}) \to 0\;.
\]
Taking into account that $u_{n}\in \nhlo$, by \eqref{fprol}, one 
deduces that $u_{n}$ strongly 
converges to zero, which is an evident contradiction with
\eqref{cesequence} implying that $(\mu_{n})$ is bounded.
As a consequence, there exists $\mu\in \R$ such that, up to a subsequence, $\mu_{n}\to \mu$. Assume, by contradiction, that $\mu\neq 0$, take $\vfi=u_{n}$ as test function in 
\eqref{equavinc} and, recalling that $u_{n}\in \nhlo$, we have
$$
\langle N'_{\Omega}(u_{n}),u_{n}\rangle\to 0.
$$
From this point we can repeat the same argument, starting from
\eqref{nprime}, to get a contradiction.
Then $\mu_{n}\to 0$ and \eqref{equavinc} yields  the conclusion. 
\edim
The next lemma studies the asymptotic behavior of a bounded
Palais-Smale sequence of $I_{\Omega}$. 
In the case $f(t)=t^{p}$, the proof
is given in \cite{bece} (see also Chapter 8 in \cite{wi}).
However, thanks to hypotheses \eqref{fprol}, \eqref{fdergrowth}
the proof can be handled as in the polynomial case arguing as in Chapter 8 in \cite{wi}. We will include some details for the sake of clearness. 

\begin{lemma}\label{splitting4}
(Splitting) 
Assume \eqref{fprol}, \eqref{fdergrowth}.
Let $(u_n)\in  H^1_0(\Omega)$ be a bounded sequence such that 
\beq\label{iposplit}
I_\Omega(u_n)\rightarrow d >0 \quad \text{and} \quad  I_\Omega'(u_n)\rightarrow 0\;\; \text{in} \;  H^{-1}(\Omega).
\eeq
Replacing $(u_n)$ by a subsequence, if necessary, there exists a solution $u_{0}$ of $\eqref{P}$, such that the following alternative holds:
\\
either $u_{n}\to u_{0}$ strongly in $\huno$ or there exist
a finite sequence 
$(u^1,u^2,...,u^k)$ in $H^{1}(\mathbb{R}^{N})$ 
solutions of \eqref{Pinf}
and $k$ sequences of points $(y^{j}_{n})\subset \mathbb{R}^N, 1\leq j \leq k$, satisfying: 
\begin{itemize}
\item[a)] $|y^{j}_{n}|\rightarrow \infty$ and $|y^{j}_{n} - y^{i}_{n}|\rightarrow \infty, i\neq j$;
\item[b)] $\|u_n -u_{0}- \displaystyle\sum_{j=1}^{k}u^j(\cdot-y^{j}_{n})\|\to 0$;
\item[c)]  $I_\Omega(u_n)\rightarrow I_\Omega(u_{0}) + 
\displaystyle\sum_{j=1}^{k}I(u^j)$.
\end{itemize}
\end{lemma}
\bdim
Since $(u_n)$ is bounded,  there exists $u_0\in  H^1_0(\Omega)$ such that, up to a subsequence, $u_n\rightharpoonup u_0$. Then, thanks to the continuity of 
$f$ and to \eqref{fdergrowth}, from \eqref{iposplit} 
it follows that 
\beq\label{equauzero}
I'_{\Omega}(u_{0})=0\qquad \text{in } H^{-1}(\Omega).
\eeq 
Let $u_n^1:=u_n-u_0 \in H^1_0(\Omega) \subset H^1(\mathbb{R}^N)$, and let us show that it results, for $n\rightarrow\infty,$
\begin{align}
\label{norm} \|u_n^1\|^2 &=\|u_n\|^2-\|u_0\|^2+o_n(1);
\\
\label{functional} I(u_n^1)&\rightarrow d-I_\Omega(u_0);
\\
\label{derivative} I'(u_n^1)&\rightarrow 0 \quad \text{in $H^{-1}(\Omega)$}.
\end{align}

The proof of \eqref{norm} is standard.
To show \eqref{functional}, 
note that the weak convergence of  $(u_n)$ to $u_0$ implies $u_n^1\rightharpoonup 0$ weakly in $H^{1}_{0}(\Omega)$.
Then,  applying  Theorem 2 in  \cite{brli} (see also Lemma 8.1 in \cite{wi}),  it follows that
\[
I(u^1_n)-I_\Omega(u_n)+I_\Omega(u_0){}=-\int_{\Omega}\left[F(u^{1}_n)-F(u_n)+F(u_0)\right]+o_{n}(1)
=o_n(1),
\]
where $o_{n}(1)$ is a quantity tending to zero as $n$ goes to
plus infinity. Then, \eqref{iposplit}   yields \eqref{functional}.
Moreover, \eqref{derivative} follows from the following facts:
observe that condition \eqref{fdergrowth}
allows us to use Theorem 2 in \cite{brli} or arguing as in
 Lemma 8.1 in \cite{wi}, to obtain that
\(f(u_{0}+u^{1}_{n})-f(u^{1}_{n})\to f(u_{0}) \) in 
\( H^{-1}(\Omega)\).
Then, exploiting \eqref{iposplit}, 
\eqref{equauzero}, it follows, for every $\vfi\in H^{1}_{0}(\Omega)$
\begin{eqnarray*}
\eps_{n}\|\vfi\|_{\Omega}&\geq &\left|\left\langle I_\Omega'(u_n),\varphi
\right\rangle \right|
=
\left|\left\langle I_\Omega'(u_0+u^1_n),\varphi\right\rangle\right|
\geq\Big| \left\langle I_\Omega'(u_0),\varphi\right\rangle 
+\left\langle I_\Omega'(u^{1}_n),\varphi\right\rangle \Big|
\\
&&-\Big|\int_{\Omega}\left[f(u_0+u^1_n)-f(u_0) -f(u^1_n)\right]\varphi \Big|
\\
&=& \left|\left\langle I'(u^1_n), \varphi\right\rangle
\right|
-\left|\int_{\Omega}\left[f(u_0+u^1_n)-f(u_0) -f(u^1_n)\right]\varphi \right|
\\
&\geq & \left|\left\langle I'(u^1_n), \varphi\right\rangle\right|
-\eps_n\|\vfi\|_{\Omega}.
\end{eqnarray*}
Hence, \eqref{derivative} follows.
Let us now consider
$$
\delta:=\limsup_{n\rightarrow\infty}\sup_{y\in
\mathbb{R}^N}\int_{B_1(y)}|u_n^1|^2.
$$ 
If $\delta=0$, it follows from Lions' Lemma  \cite{lions} that
\begin{equation}\label{ss9}
u_n^1\rightarrow 0 \quad \mbox{in} \ L^{p}(\mathbb{R}^N), \; \mbox{for} \ \mbox{any $p$ }\,:\,  2<p<2^*.
\end{equation}
On the other hand, as $(u^{1}_{n})$
 is bounded in $H^{1}_{0}(\Omega)$,  we can use 
\eqref{fdergrowth} and \eqref{ss9}
to deduce that $u^{1}_{n}\to 0$ strongly in $H^{1}_{0}(\Omega)$
and since $u^{1}_{n}\equiv 0$ outside $\Omega$ it results
$u^{1}_{n}\to 0$ strongly in $H^{1}(\R^{N})$.
Then the first alternative in the statement of the Lemma holds. 
While, if $\delta>0$, there exists a sequence $(y_n^{1})\subset
\mathbb{R}^N$ such that
\begin{equation}\label{ss11}
\int_{B_1(y^{1}_n)}|u_n^1|^2>\frac{\delta}{2}.
\end{equation}
Define a new sequence $(v_n^1)\subset H^1(\mathbb{R}^N)$ by
$v_n^1:=u_n^1(\cdot+y^{1}_n).$ Since $(u_n^1)$ is bounded in 
$H^{1}(\R^{N})$, then $(v_n^1)$ is also bounded and there 
exists $u^{1}\in H^1(\mathbb{R}^N)$, such that
$v_n^1\rightharpoonup u^1$ in $H^1(\mathbb{R}^N)$
and $v_n^1\rightarrow u^1$ almost everywhere in $ \mathbb{R}^N$; then
performing a change of variable  in  \eqref{ss11}
and applying Fatou Lemma one obtains that 
\[
\int_{B_1(0)}|v_n^1|^2>\frac{\delta}{2}\;,
\qquad
\int_{B_1(0)}|u^1|^2\geq\frac{\delta}{2}\;,
\]
so that $u^1\not\equiv 0$. Moreover, since
$u_n^1\rightharpoonup 0$ in $H^1(\mathbb{R}^N)$, it follows that 
up to a subsequence we can assume that
$|y^{1}_n|\rightarrow\infty$.
As shown in \cite{bece} (see also Proposition 3.1 in \cite{cemopa}
or Chapter 8 in \cite{wi}),  it results that $I'(u^1)=0$ 
in $H^{-1}(\R^{N})$. 
Indeed, take  $\phi \in C^\infty_c(\R^{N})$ and observe that,
since $|y_{n}^{1}|\to+\infty$ we can find $n_{0}$ such that
$\phi_{n}:=\phi(x-y_{n}^{1})\in C^{\infty}_{c}(\Omega)$ 
 for every $n\geq n_{0}$;
moreover,  
$\|\phi_{n}\|_{\Omega}\leq \|\phi\|$.
As a consequence, \eqref{derivative} yields
\[
\begin{split}
\left|\langle I'(u^1),\phi\rangle\right| &=
\left|\langle I'(v^1_n),\phi\rangle\right|+o_{n}(1)
=  \left|\langle I'(u^1_n),\phi_{n}\rangle\right|+o_{n}(1)
=o_{n}(1),
\end{split}
\]
so that $u^{1}$ is a solution of \eqref{Pinf}.
Defining  $\dys  u_n^2(x):=u_n^1(x)-u^1(x-y_n^1),$ it results that
 $u_{n}^{2}(\cdot+y^{1}_{n})=v^{1}_{n}-u^{1}$ and
repeating for $u^{2}_{n}$ the argument done for $u^{1}_{n}$
we deduce that
$$
\|u^{2}_{n}\|^{2}=\|u_n^1\|^2-\|u^1\|^2+o_n(1),
\quad  I(u^{2}_{n})\to  d-I_{\Omega}(u_{0})-I(u^{1}),
\quad
I'(u^{2}_{n})\to 0 \;\text{in $H^{-1}(\Omega)$},
$$
so that we go on repeating the argument obtaining $(y^{2}_{n})$
satisfying conclusion a) and  $u^{2}$,
another solution of \eqref{Pinf}.
From now on
we proceed by iteration. Note that if $u$ is 
a nontrivial critical point of $I$ and $w$ is the solution  of 
minimum action of  $(\ref{Pinf})$, then 
\begin{equation}\label{i112}
I(u)\geq I(w)>0.
\end{equation}
As a consequence, passing from the step $l$ to the step $l+1$ 
the action level decreases 
because  in the asymptotic information on the 
functional $I(u^{l+1}_{n})$ it appears $-I(u^{l+1})$.
Then, taking into account \eqref{i112}, 
the sum must have a finite number of terms, so that,
the iteration must be finite and terminate at some index 
$k\in \mathbb{N}$, yielding also conclusions b) and c).  
\edim
The following  result is a direct  consequence of the Splitting Lemma.
\begin{corollary}\label{ceramirange}
Assume \eqref{fprol}, \eqref{fdergrowth}, \eqref{fmono},  
\eqref{finfty}, \eqref{NQ}, \eqref{iponh} and \eqref{unique}.
Let $(u_{n})$ be a  sequence sa\-tisfying \eqref{cesequence}
and \eqref{ipo} with $d\in (m,2m)$. 
Then $(u_{n})$ admits a strongly convergent subsequence.
\end{corollary}
\bdim
Applying Lemma \ref{ceramibounded} we deduce that $(u_{n})$ is bounded, moreover Lemma \ref{nuovo} allows us to apply 
Lemma \ref{splitting4} to obtain that there exists 
$u_{0}$ solution of \eqref{P} such that
conclusion {\it b)} holds. Moreover,
conclusion {\it c)} and Lemma \ref{equal} yield
\[
2m>d=\sum_{j=0}^{k}I(u^{j})\geq
\begin{cases} 
km &\text{if $u_{0}\equiv 0$,}
\\
m_{\Omega}+km= (k+1)m &\text{if $u_{0}\not\equiv 0$.}
\end{cases}
\]
Then, in both cases, $k<2$, i.e. $k=0$
or $k=1$.
If $k=1$ and $u_{0} \equiv 0$, it follows that  $I_{\Omega}(u_n)\rightarrow I(u_1)=d$ and from hypothesis \eqref{fprol}, $u_{1}$ 
is positive, so that, condition \eqref{unique} yields
$u_{1}=w$ and $d=I(u_{1})=m$, which contradicts the hypothesis. 
Also note that the hypothesis $d<2m$ implies that
it is not possible that $k=1$ and $u_{0}\not\equiv 0$.
Then, the only possibility is that $k=0$,  that is
$d=I(u_{0})=I_{\Omega}(u_{0})$ and
 Lemma \ref{splitting4} implies
$u_{n}$ strongly converges to $u_{0}$.
\edim
\section{The Linking Argument and Proof of the Main Results}
\label{final}
We will prove our existence results by applying the
 Linking theorem on the manifold $\nhlo$ 
(see Theorem 8.22 in \cite{amma} or Theorem 8.4 in \cite{struwe}
joint with Lemma 5.14 and 5.15 in \cite{wi}).
This argument has already been used in \cite{amceru}
to prove existence results for linearly coupled semi-linear
non-autonomous equations in $\R^{N}$.
In order to define the linking sets, we will make use of 
the properties of a  barycenter function,
already adopted in \cite{bece} and then often used
when building solution at higher action level than 
the least one (see for example  \cite{amceru, cepa95}
or \cite{bawe} for an interesting generalization). 
Here, we will follow the notation in \cite{cepa95} (see 
also \cite{amceru}).
For every $u\in H^{1}(\R^{N})\setminus \{0\}$, the 
following maps are well defined
\[
\mu(u)(x):=\frac{1}{|B_{1}(x)|}\int_{B_{1}(x)}|u(y)|dy,\quad
\text{$\mu(u)\in L^{\infty}\cap C^{0}(0,+\infty)$ },
\]
$$
\hat{u}(x):=\left[\mu(u)(x)-\frac{\|\mu(u)\|_{\infty}}2\right]^{+},
\quad \hat{u}\in C_{0}(\R^{N}).
$$
Then, the barycenter  of a function $u\in H^{1}(\R^{N})\setminus \{0\}$  defined by
$$
\beta(u)=\frac{1}{\|\hat{u}\|_{1}}\intr x\hat{u}(x)dx
$$
 is a continuous function enjoying the following properties 
\begin{align}
\label{btras} \beta(u(\cdot-y))&=\beta(u)+y\quad \forall\,y\in \R^{N},
\\
\label{bmolt}\beta(Tu)&=\beta(u)\quad \forall\,T>0.
\end{align}
One of the linking set is the following subset of
${\mathcal N}_{\Omega}$
\beq\label{defS}
S:= \left\lbrace u\in H^{1}_0(\Omega),\,:\; u\in \nhlo, \; \beta(u)=0\right\rbrace\;.
\eeq
\bos\label{remb}
Notice that $S\neq \emptyset$. 
Indeed, first note that, 
from the properties of the function $w$ we can find 
$r_{0}>0$  such that 
$\{x\in \R^{N}\,:\, \mu(w)>\|\mu(w)\|_{\infty}/2\}= B_{r_{0}}(0)$.
Then, take $\theta\in \R^{N}$ with $|\theta|=1$, $R>4\rho r_{0}$
and the functions $\eta_{1}$ and $\eta_{2}$ defined as  
$\eta_{1}(x)=1-\xi(4\rho(x-R\theta)/R)$
and $\eta_{2}(x)=1-\xi(4\rho(x+R\theta)/R)$. Choose the function
$z_{R}=w^{R\theta}\eta_{1}+w^{-R\theta}\eta_{2}$.
The properties of the function $\eta_{1}$ imply that
 $\eta_{1}w^{R\theta}<w^{R
\theta}$ so that $\mu(\eta_{1}w^{R\theta})
<\mu(w^{R\theta})$ and it results
\[
\begin{split}
\sup_{\R^{N}\setminus B_{R/4\rho}(R\theta)}\mu(\eta_{1}w^{R\theta})
<\sup_{B_{R/4\rho}(R\theta)}\mu(\eta_{1}w^{R\theta})
=
\sup_{ B_{R/4\rho}(R\theta)}\mu(w^{R\theta})=\mu(w^{R\theta})(R\theta),
\end{split}
\]
where we have taken into account that $\mu(0)>\mu(x)$
for every $x\in \R^{N}\setminus\{0\}$.
As a consequence of the previous facts, 
\[
\left\{x\in \R^{N}\,:\, \mu(\eta_{1}w^{R\theta})>\frac{\|\mu(\eta_{1}w^{R\theta})\|_{\infty}}{2}\right\}\subset B_{R/4\rho}(R\theta).
\]
Arguing analogously for $\eta_{2}w^{-R\theta}$, one obtains
 $\widehat{\eta_{1}w}^{R\theta}=\widehat{w}^{R\theta}$
and $\widehat{\eta_{2}w}^{-R\theta}=\widehat{w}^{-R\theta}$.
Moreover, observe that $\GO(z_{R})=\GO(w^{R\theta})
+\GO(w^{-R\theta})+ o(R)=2\GO(w)+o(R)$, 
where $o(R)$ is a quantity tending to zero as $R$ tends to infinity.
Then $z_{R}\in A_{\Omega}$ 
for $R$ sufficiently large.
Therefore, $v=T_{\Omega}(z_{R})z_{R}\in {\mathcal N}_{\Omega}$ 
and from \eqref{bmolt} it results $\beta(v)=\beta(z_{R})$.
In addition taking into account that $B_{R/2}(R\theta)\cap
B_{R/2}(-R\theta)=\emptyset$, we obtain
$\hat{z}_{R}=\widehat{\eta_{1}w}^{R\theta}+\widehat{\eta_{2}w}^{-R\theta} $. Then it results
\[
\begin{split}
\beta(v)=\beta(z_{R})
=
\beta(\eta_{1}w^{R\theta})+\beta(\eta_{2}w^{-R\theta})=\beta(w^{R\theta})+\beta(w^{-R\theta})=
0,
\end{split}
\]
showing that $z_{R}\in S$.
\eos
\begin{lemma} \label{bbound}
Assume \eqref{fprol}, \eqref{fdergrowth}, \eqref{fmono},
\eqref{finfty}, \eqref{NQ}, \eqref{iponh}.
Let $b = \displaystyle \inf_{S}I_\Omega(u)$, then $b > m_{\Omega}$.
\end{lemma}
\bdim
It is clear that $b \geq m_{\Omega}$.
To prove the strict inequality we shall argue by contradiction. 
Suppose $b = m_{\Omega}$, then there exists a sequence $(v_n) \subset \nhlo$
such that $\beta(v_n)=0$, $\langle I'_\Omega(v_n), v_n\rangle =0$ and $I_\Omega (v_n) \to m_{\Omega}$; 
moreover, $(v_n)$ is not relatively compact because $m_{\Omega}$ is not attained.
By Ekeland Variational Principle (Theorem 8.5 in \cite{wi}) applied to  the closed manifold $\nhlo$, 
there exists another sequence
$(\tilde v_n) \subset \nhlo$ such that:
\beq\label{vinc}
I_{\Omega}(\tilde v_n) \to m_{\Omega}\;,
\quad
 \nabla_{\nhlo} I_{\Omega}(\tilde v_n)\to 0 \;,
\quad
\| \tilde v_n- v_n \|_{\Omega} \to 0 \;.
\eeq
From Lemma \ref{ceramibounded} $(v_{n})$ and $(\tilde{v}_{n})$
are bounded.  Moreover, Lemma \ref{nuovo} implies 
that  
\beq\label{eq:free}
I_{\Omega}^{\prime}(\tilde v_n) \to 0.
\eeq
By exploiting  hypothesis \eqref{fdergrowth} we deduce that
that $I''_{\Omega}$ maps bounded sets of $H^{1}_{0}(\Omega)$
in bounded sets, then
the mean value Theorem implies that the following
inequality holds for every $\phi\in H^{1}_{0}(\Omega)$, 
$$
\left| 
\langle I_{\Omega}'(\tilde{v}_{n}) - I_{\Omega}'(v_{n}), 
\phi\rangle 
\right| 
\leq K  \| \tilde{v}_{n}
 -v_{n}\|_{\Omega}  \|\phi\|_{\Omega} \,.
$$
Taking the supremum on $\phi$ and using \eqref{eq:free},
it follows that also $(v_{n})$ is 
a   Palais-Smale sequence.
Therefore, from Lemma \ref{splitting4} and since
$m_{\Omega}$  is not attained, we deduce that
conclusions {\it a)-c)} hold.
In particular, conclusion {\it c)}, \eqref{vinc} 
and Lemma \ref{equal} imply 
$$
m=I_{\Omega}(u_{0})+\sum_{j=1}^{k}I(u^{j})\geq 
I_{\Omega}(u_{0})+ km\geq km,
$$
where the last inequality is implied by the fact that, either
$u_{0}\equiv0$ so that $\io(u_{0})=0$, or $u_{0}\in \nhlo$ so that $I(u_{0})>m_{\Omega}>0$.
Then,  $k$ has to be equal to one, $u_{0}$ has to be
trivial and $u^{1}=w$. 
This and conclusion {\it b)} of Lemma \ref{splitting4} yield
$v_n(\cdot)-w(\cdot -y_n)\rightarrow 0$, strongly in 
$H^{1}(\R^{N})$, where $y_n \in \mathbb{R}^N$, 
$|y_n| \to \infty$.
Calculating the barycenter function of $v_{n}(x)$ and
$w(x-y_{n})$, we have,
as $\beta$ is a continuous function,  $w$ 
is radially symmetric,  and using \eqref{btras}
$$
0=\beta( v_n) =\beta( w(\cdot -y_n))+o(1)
=\beta(w)+y_{n}+o(1)=
y_{n}+o(1) \;
$$
where $o(1)$ a quantity tending to zero as $n$ goes to plus
infinity.
This  immediately gives a contradiction as
$|y_{n}|\to +\infty$.
\edim
In order to define the other linking set, 
we argue as in \cite{cepa95}
(see also Section 7 in \cite{amceru})
and we take $x_{0}$ and $y$ as in \eqref{def:ur} and 
$R>\max\{R_{0},R_{1}\}$  where $R_{0}$ and 
$R_{1}$ are  introduced in Lemmas \ref{tprop} and 
\ref{le:uniform}. Let us define the function 
$$
z: \partial B_2(x_{0}) \times [0,1] \mapsto \ov{B_{2}(x_{0})},\quad
\text{ by }\quad
z(y,t):= t y + (1-t) x_{0}.
$$ 
Notice that 
$z$ is a homeomorphism from $\partial B_2(x_{0}) \times (0,1]$ to 
$ \overline{B_{2}(x_{0})}\setminus \{x_{0}\}$, 
so that for every point $\tilde{z}\in \ov{B_{2}(x_{0})}$
there is a unique pair $(y,t)\in \partial B_2(x_{0}) \times (0,1] $,
such that $\tilde{z}=z(y,t) \not = x_{0}$. 
Therefore, we can define the  operator
$\Psi_{R}:\overline{B_2(x_{0})} \mapsto \nhlo$ by the values that it 
takes on $\partial B_{2}(x_{0})\times [0,1]$ as follows
\beq\label{def:psi}
\Psi_{R}[z]=\begin{cases}
\Psi_{R}[z(y,t)]=\Psi_{R}[y,t]=
\Pi_{\nhlo}(U^{R}_{t})\; &\text{if $z\neq x_{0}$},
\medskip\\
 \Psi_{R}[y,0]=\Pi_{\nhlo}(\Phi^{Rx_{0} }) &\text{if $z=x_{0}$},
\end{cases}
\eeq
where $U^{R}_{t}$ is given in \eqref{def:ur} and $\Pi_{\nhlo}$
is defined in \eqref{defpi}.
Notice that  $\Psi_{R}$
is well defined thanks to Lemma  \ref{le:uniform}.
In order to apply the Linking Theorem on $\nhlo$
we recall that $S$ is defined in \eqref{defS} and we set  
\begin{align}
\label{defq} Q:=\Psi_{R}(\overline{B_2(x_{0})})&
\\
\label{defh}{\mathcal H}: = \left\lbrace h\in C^0(Q,{\mathcal N}_{\Omega}), \;\, 
h\vert_{\partial Q} = \text{id}\right\rbrace, &\quad 
d: = \displaystyle\inf_{h\in \MH}\sup_{u\in Q}\io(h(u)).
\end{align}
The linking geometrical structure is showed in the next lemma.
\begin{lemma}\label{geometry}
Assume conditions \eqref{fprol}, \eqref{fdergrowth},
\eqref{fmono}, \eqref{finfty},
\eqref{NQ}, \eqref{iponh}.
Let $Q$ and ${\mathcal H}$ be defined in \eqref{defq} \eqref{defh}.
Then, for $R$ sufficiently large, it results 
\begin{align}
\label{cap} \partial Q\cap S &=\emptyset,
\\
\label{hmap}h(Q)\cap S &\neq \emptyset \,\qquad \forall h\in \MH,
\\
\label{supinf}
\sup_{\partial Q}I_{\Omega}(u) &<\inf_{S}I_{\Omega}(u).
\end{align}
\end{lemma}
\bdim
First of all, notice that  $\Psi_{R}:\ov{B_{2}(x_{0})}\mapsto
\Psi_{R}\left(\ov{B_{2}(x_{0})}\right)$
is a continuous bijection defined on a compact set.
Indeed, let us first show that $\Psi_{R}$ is injective 
in $\ov{B_{2}(x_{0})}\setminus \{x_{0}\}$. 
In order to do this, let us consider $z_{1},\,z_{2}\in
\ov{B_{2}(x_{0})}\setminus \{x_{0}\}$ such 
that $\Psi_{R}(z_{1})=\Psi_{R}(z_{2})$. Since
$z(y,t)$ is injective, this amounts to consider
$(y_{1},t_{1}),\, (y_{2},t_{2})\in \partial B_{2}(x_{0})
\times (0,1]$  such that $\Psi_{R}[z(y_{1},t_{1})]=\Psi_{R}[z(y_{2},t_{2})]$.
Taking into account Remark \eqref{defpi}, and using the notation
\[
T_{i}(R)=T_{\Omega}(t_{i} \Phi^{Ry_{i}}+(1-t_{i})\Phi^{Rx_{0}}),
\]
this is equivalent to have
\[
T_{1}(R)\left(t_{1} \Phi^{Ry_{1}}+(1-t_{1})\Phi^{Rx_{0}}\right)=
T_{2}(R)\left(t_{2} \Phi^{Ry_{2}}+(1-t_{2})\Phi^{Rx_{0}}\right).
\]
Suppose that $|y_{1}-y_{2}|=a>0$. Since
$|y_{i}-x_{0}|=2$ for $i=1,2$,  
\eqref{behavior:w} implies that 
$w(R(y_{2}-y_{1}))=o(R)$
and $w(R(y_{i}-x_{0}))=o(R)$,
where $o(R)$ is a quantity tending to zero as $R$ tends to infinity.
Choosing 
$x=Ry_{1}$ with $R> 2\rho$,
and recalling \eqref{phi}, one obtains
\(
T_{1}(R)\left(t_{1}w(0)+(1-t_{1})o(R)\right)=T_{2}(R)o(R).
\)
Then, \eqref{uniform} implies that
$ T_{1}(R)t_{1}w(0)= o(R).$
Since $w(0)>0$ and $t_{1}\in (0,1]$, this implies that
\(T_{1}(R)=o(R).\)
Arguing as in the proof of Lemma \ref{le:uniform}, one can reach a contradiction,
so that $|y_{1}-y_{2}|=0$, i.e. $y_{1}=y_{2}$.
In order to prove that $t_{1}=t_{2}$ let us  choose  $x=Ry=Ry_{1}=Ry_{2}$, take into account \eqref{defxi} and that $R>2\rho$ to  obtain
\[
T_{1}(R)\left[t_{1}w(0)+(1-t_{1})w(2R)\right]=
T_{2}(R)\left(t_{2}w(0)+(1-t_{2})w(2R)\right).
\]
On the other hand, choosing $x=Rx_{0}$ we get
\[
T_{1}(R)\left[ t_{1} w(2R)+(1-t_{1})w(0)\right]=
T_{2}(R)\left[t_{2} w(2R)+(1-t_{2})w(0)\right].
\]
Then
\[
\dfrac{t_{1}w(0)+(1-t_{1})w(2R)}{ t_{1} w(2R)+(1-t_{1})w(0)}=
\dfrac{t_{2}w(0)+(1-t_{2})w(2R)}{t_{2} w(2R)+(1-t_{2})w(0)}.
\]
Since the function $h(t)=at+b(1-t)/[bt+a(1-t)]$ is injective,
the above equality implies that $t_{1}=t_{2}$. 
Then, we have shown that $\Psi $ is injective in $\overline{B_{2}(x_{0})}\setminus\{x_{0}\}$. Now, take $z\in \overline{B_{2}(x_{0})}$ such that  $\Psi_{R}(z)=\Psi_{R}(x_{0})$, that is 
\[
T_{\Omega}\left(t \Phi^{Ry}+(1-t)\Phi^{Rx_{0}}\right)\left[t\Phi^{Ry}+(1-t)\Phi^{Rx_{0}}\right]=
T_{\Omega}\left(\Phi^{Rx_{0}}\right) \Phi^{Rx_{0}}.
\]
Choosing $x=Ry$ and arguing as before, one obtains that $z=x_{0}$
proving that $\Psi_{R}$ is injective in $\overline{B_{2}(x_{0})}$. As a 
consequence, $\Psi_{R}$ is an homeomorphism and $\partial Q=\Psi_{R} (\partial B_{2}(x_{0}))$.\\
In order to show that \eqref{cap} holds, we observe that,
since $w$ is radially symmetric, positive and decreasing in 
$(0,+\infty)$, also $\mu(w)$ is decreasing with respect to $|x|$;
Moreover, as proved in Theorem 2.1 in \cite{bawe}
$\mu(w)\to 0$ as $|x|\to +\infty$, then,
arguing as in Remark \ref{remb}, we obtain that
there exists a unique $r_{0}>0$ such that for every 
$|x|=r_{0}$, $\mu(w)(x)=\|\mu(w)\|_{\infty}/2$
and by \eqref{btras}
the set
\[
E(w):=\left\{x\in \R^{N},
\mu(w)(x)\geq \frac{\|\mu(w)\|_{\infty}}{2}\right\}
\]is such that
\beq\label{defbr}
E(w)=B_{r_{0}}(0) \Rightarrow E(w(\cdot-Ry))=B_{r_{0}}(Ry), 
\eeq
for every $ R\in \R^{+}$.
Let us fix $R$ such that $ R>2\rho+1+r_{0}$ and,
as $y\in\partial B_{2}(x_{0})$, it results that  
$|x|>2\rho+r_{0}$ for every
$x\in B_{1}(Ry)$. Then, exploiting \eqref{defxi} we obtain
\begin{align*}
\mu(\Phi^{Ry})(Ry)&=\frac{1}{|B_{1}(Ry)|}
\int_{B_{1}(Ry)}\xi(x)w(x-Ry)dx=
\frac{1}{|B_{1}(Ry)|}\int_{B_{1}(Ry)}w(x-Ry)dx
\\
&=\frac{1}{|B_{1}(0)|}\int_{B_{1}(0)}w(\sigma)d\sigma=\mu(w)(0)
=\|\mu(w)\|_{\infty}.
\end{align*}
Since, $\|\xi\|_{\infty}\leq 1$, it results
\beq\label{uffa}
|\mu(\Phi^{Ry})(x)|\leq |\mu(w)(x-Ry)|\leq \|\mu(w)\|_{\infty},
\eeq
showing that 
\beq\label{eq:uffa2}
\|\mu(\Phi^{Ry})\|_{\infty}=|\mu(\Phi^{Ry})(Ry)|=
\|\mu(w)\|_{\infty}.
\eeq 
In addition,  for every $x\in B_{r_{0}}(Ry)$,
any $z\in B_{1}(x)$ satisfies
$|z|>2\rho,$ showing that
$B_{1}(x)\subseteq \R^{N}\setminus B_{2\rho}(0)$, 
and using again \eqref{defxi}, we have
\[
\begin{split}
\mu(\Phi^{Ry})(x)
&=\frac{1}{|B_{1}(x)|}
\int_{B_{1}(x)}\xi(z)w(z-Ry)dz=\frac{1}{|B_{1}(x)|}
\int_{B_{1}(x)}w(z-Ry)dz
\\
&=\mu(w)(x-Ry).
\end{split}
\]
Recalling \eqref{defbr} we have that, for every
$x\in B_{r_{0}}(Ry)$, 
$\mu(w)(x-Ry)>\|\mu(w)\|_{\infty}/2$, yielding
$$
\mu(\Phi^{Ry})(x)>\frac12 \|\mu(w)\|_{\infty}\qquad \text{for
every $x\in B_{r_{0}}(Ry)$},
$$
so that  $\widehat{\Phi}^{Ry}\neq 0$ if $x\in B_{r_{0}}(Ry)$. 
If $x\not\in B_{r_{0}}(Ry)$ recalling  \eqref{defbr}, \eqref{uffa} 
and \eqref{eq:uffa2}, it results
$$
\mu(\Phi^{Ry})\leq\mu(w(\cdot-Ry))<\frac12 \|\mu(w(\cdot-Ry))\|_{\infty}=\frac12 \|\mu(\Phi^{Ry})\|_{\infty}.
$$
Therefore, $\widehat{\Phi}^{Ry}\neq 0$ if and only if 
$x\in B_{r_{0}}(Ry)$,
but, in this set $\xi\equiv 1$, so
that $\widehat{\Phi}^{Ry}=\widehat{w}(\cdot-Ry)$ and hence
\beq\label{betauffa}
\beta(\Phi^{Ry})=\beta(w(\cdot-Ry))=Ry,
\eeq
showing \eqref{cap}.
In order to prove \eqref{hmap}, for every $h\in {\mathcal H}$, 
let us take $V:\ov{B_{2}(x_{0})}\mapsto \R^{N}$ given by 
$V(z)=\left(\beta \circ h\circ \Psi_{R}\right) (z)$.
If $z\in \partial B_{2}(x_{0})$, 
$z=z(y,1)$, then, it results
$\Psi_{R}(z)\in \partial Q$, so that 
$h(\Psi_{R}(z))=\Psi_{R}[y,1]$, 
then \eqref{bmolt} and \eqref{betauffa} yield
\[
V(z)=\beta(\Psi_{R}[y,1])=\beta(\Psi_{R}(\Phi^{Ry}))=
\beta(T_{\Omega}(\Phi^{Ry})\Phi^{Ry})=
\beta(\Phi^{Ry})=Ry.
\]
Then, by Brouwer Theorem, there exists $z_{0}\in B_{2}(x_{0})$ such that  $V(z_{0})=0$,  i.e. $h(\Psi_{R}(z_{0}))\in S$, and, as by definition $\Psi_{R}(z_{0})\in Q$, 
\eqref{hmap} follows.
\\
In order to show \eqref{supinf}, let us first observe that
$\partial Q=\Psi_{R}(\partial B_{2}(x_{0}))=\Psi_{R}[y,1]$;
moreover,  
$$
I_{\Omega}(\Psi_{R}[y,1])=I_{\Omega}(\Pi_{\nhlo}(\Phi^{Ry}))=I_{\Omega}(T_{\Omega}(\Phi^{Ry})\Phi^{Ry}).
$$ 
From \eqref{limyI} and 
\eqref{limtn} it follows that 
$I_{\Omega}(\Psi_{R}[y,1])\to m$ as $R\to+\infty$; 
then, Lemma  \ref{equal} and Lemma \ref{bbound} imply 
that, fixing $R$ sufficiently large, \eqref{supinf} holds.
\edim
\ble\label{prop:chiave}
Assume \eqref{fprol}, \eqref{fdergrowth}, \eqref{fmono}, \eqref{finfty}.
There exists $R_{2}>\max\{R_{0},R_{1}\}$  such that 
\[
\max_{Q}\io(u)<2m_{\Omega},
\quad \text{for every $R\geq R_{2}$}.
\]
\ele
\bdim
In the proof we will  use the notation, introduced in Lemma \ref{ack}
\[
T^{R}_{t}=T_{\Omega}(U^{R}_{t}),
\]
where $T(u)$ is defined in Lemma \ref{tprop}.
Recalling \eqref{defq}, \eqref{def:psi}, and taking into account
\eqref{def:ur},  it is sufficient to  show that there exists $R_{2}$ 
sufficiently large, such that
\beq\label{eq:dim}
\io(\Pi_{\nhlo}(U^{R}_{t}))=\io(T^{R}_{t}U^{R}_{t})<2m_{\Omega}, \quad \forall\, R\geq R_{2}, \, \forall\,(y,t)\in \partial B_{2}(x_{0})\times [0,1].
\eeq
It results
\begin{align*}
\io(T^{R}_{t}U^{R}_{t})
=&
\frac{(T^{R}_{t}t)^{2}}{2}\|\firy\|_{\Omega}^{2}
+\frac{(T^{R}_{t}(1-t))^{2}}{2}\| \firxz\|_{\Omega}^{2}
+(T^{R}_{t})^{2}t(1-t)\into \nabla \firy\nabla \firxz 
\\
&+(T^{R}_{t})^{2}t(1-t)\la\into \firy\firxz
-\into F(T^{R}_{t}t\firy+T^{R}_{t}(1-t)\firxz).
\end{align*}
By adding and subtracting $I(T^{R}_{t}tw^{Ry})
+I(T^{R}_{t}(1-t)w^{Rx_{0}})$ and 
taking into account
\eqref{limy} and Lemma \ref{le:uniform}
one obtains
\[
\begin{split}
\io(T^{R}_{t}U^{R}_{t})=
& 
I(T^{R}_{t}tw^{Ry})+I(T^{R}_{t}(1-t)w^{Rx_{0}})
+
\intr \left[F(T^{R}_{t}tw^{Ry})-F(T^{R}_{t}t\firy)\right]
\\
&-\intr \left[F(T^{R}_{t}U^{R}_{t})
-F(T^{R}_{t}t\firy)- F(T^{R}_{t}(1-t)\firxz)\right]
\\
&
+(T^{R}_{t})^{2}t(1-t)\into \nabla \firy\nabla \firxz
\la \firy\firxz +o(\eps_{R})
\\
&+\intr \left[
F(T^{R}_{t}(1-t)w^{Rx_{0}})-F(T^{R}_{t}(1-t)\firxz)
\right].
\end{split}
\]
Moreover, using \eqref{Fprop} and \eqref{uniform}, one obtains
\begin{align*}
\io(T^{R}_{t}U^{R}_{t})
&\leq I(T^{R}_{t}tw^{Ry})+I(T^{R}_{t}(1-t)w^{Rx_{0}}) +o(\eps_{R})
\\
&+(T^{R}_{t})^{2}t(1-t)\into \nabla \firy\nabla \firxz +\la\firy\firxz
\\
&-\intr\left[
F(T^{R}_{t}U^{R}_{t})-F(T^{R}_{t}t\firy)-F(T^{R}_{t}(1-t)\firxz)\right].
\end{align*}
Applying Lemmas \ref{scalprod} and \ref{ack} we deduce
\[
\begin{split}
&\io(T^{R}_{t}U^{R}_{t})
\leq 
I(T^{R}_{t}tw^{Ry})+I(T^{R}_{t}(1-t)w^{Rx_{0}})+o(\eps_{R})
\\
+
&(T^{R}_{t})^{2}t(1-t)\left[\frac12\intr f( w^{Ry}) w^{Rx_{0}}+\frac12\intr f( w^{Rx_{0}}) w^{Ry}
\right]
\\
-\intr &\left[f(T^{R}_{t}t\firy)T^{R}_{t}(1-t)\firxz+
f(T^{R}_{t}(1-t)\firxz)T^{R}_{t}t\firy\right]. 
\end{split}
\]
Applying Lemmas \ref{le:uniform} and \ref{tvm}
(see \eqref{stima:f}) in the last term one gets
\[
\begin{split}
\io(T^{R}_{t}U^{R}_{t})
\leq & 
I(T^{R}_{t}tw^{Ry})+I(T^{R}_{t}(1-t)w^{Rx_{0}})+o(\eps_{R})
\\
+
&(T^{R}_{t})^{2}t(1-t)\left[\frac12\intr f( w^{Ry}) w^{Rx_{0}}+\frac12\intr f( w^{Rx_{0}}) w^{Ry}
\right]
\\
-&
\intr \left[f(T^{R}_{t}tw^{Ry})T^{R}_{t}(1-t)w^{Rx_{0}}+f(T^{R}_{t}(1-t)w^{Rx_{0}})T^{R}_{t}tw^{Ry}\right].
\end{split}
\]
We use \eqref{inttvm} with  the choices
$\tau =T^{R}_{t}t$, $\tau=T^{R}_{t}(1-t)$
and use  \eqref{uniform}   to get
that there exists a positive constant $C_{1}$ such that
\[
\begin{split}
\io(T^{R}_{t}U^{R}_{t})
\leq & 
I(T^{R}_{t}tw^{Ry})+I(T^{R}_{t}(1-t)w^{Rx_{0}})+o(\eps_{R})
\\
-&\frac12
\intr \left[f(T^{R}_{t}tw^{Ry})T^{R}_{t}(1-t)w^{Rx_{0}}+f(T^{R}_{t}(1-t)w^{Rx_{0}})T^{R}_{t}tw^{Ry}\right]
\\
&+C_{1}\left[|tT^{R}_{t}-1|
+|(1-t)T^{R}_{t}-1|\right]\eps_{R}.
\end{split}
\]
Choosing $(\tau_{1},\tau_{2})=(T^{R}_{t}(1-t),T^{R}_{t}t)$
 in the first inequality in Lemma
\ref{le:lower},  
and $(\tau_{1},\tau_{2})=(T^{R}_{t}t,T^{R}_{t}(1-t))$ in the second one,  we  obtain (thanks to \eqref{uniform})
 \beq\label{ineq7}
\begin{split}
\io(T^{R}_{t}U^{R}_{t})
\leq 
I(T^{R}_{t}tw^{Ry})+I(T^{R}_{t}(1-t)w^{Rx_{0}})&+o(\eps_{R})
\\
-\eps_{R}
C_{1}\left[|tT^{R}_{t}-1|+|(1-t)T^{R}_{t}-1|\right]
&-\eps_{R}C_{2}\min\{tT^{R}_{t}, (1-t)T^{R}_{t}\},
\end{split}
\eeq
with $C_{2}$ a positive constant.
In view of  \eqref{lim:12}, there exists $\delta>0$ such that
for every $ t\in \left(\frac12-\delta,\frac12+\delta\right),$
it results 
\beq\label{eq:cfinal}
C_{1}\left[|tT^{R}_{t}-1|+|(1-t)T^{R}_{t}-1|\right]
+C_{2}\min\Big\{tT^{R}_{t}, (1-t)T^{R}_{t}\Big\}
\geq \frac{C_{2}}2.
\eeq
Moreover, the real valued function  $h_{v}$ defined as
\beq\label{def:finalh}
h_{v}(p):= \frac{p^{2}}{2}f(v)v-F(pv),\quad I(pv)=\intr h_{v}(p).
\eeq
achieves the unique  maximum value for $p=1$ (see
for more details Lemma \ref{ceramibounded} or
 Proposition  3.20 in \cite{mamope}) then, as $w\in {\mathcal N}$, 
 it results
\[
I(T^{R}_{t}tw^{Ry})
\leq \intr \frac12 f(w)w- F(w)=I(w)=m.
\]
The same conclusion holds for the term 
$I(T^{R}_{t}(1-t)w^{Rx_{0}}))$, so that
from \eqref{eq:cfinal} it follows
\begin{align}\label{ineq8}
\io(T^{R}_{t}U^{R}_{t})
\leq & 
2m-\frac{C_{2}}2\eps_{R}+o(\eps_{R})
\end{align}
which, recalling Lemma \ref{equal}, shows the conclusion if $t\in (\frac12-\delta,\frac12+\delta)$.
\\
In the other case, namely when $|t-1/2|>\delta$, we can fix a positive 
constant $\sigma_{0}$ such that 
\beq\label{tfinal}
|tT^{R}_{t}-1|> \sigma_{0}\quad \text{or}\quad
|(1-t)T^{R}_{t}-1|> \sigma_{0}.
\eeq 
Indeed, if \eqref{tfinal} were not true, we could find positive sequences
$(\sigma_{n}),\,(t_{n}),\,(T_{n})$  such that, up to  subsequences,
$\sigma_{n}\to 0$, $t_{n}\to t_{0}\in [0,1]$,
$T_{n}\to T_{0}\in [0,L]$ (thanks to \eqref{uniform}), and
$$
|t_{n}T_{n}-1|\leq \sigma_{n}\quad \text{and}\quad
|(1-t_{n})T_{n}-1|\leq \sigma_{n}.
$$
Since $\sigma_{n}\to0$, it results $(1-t_{0})T_{0}=1$
and $t_{0}T_{0}=1$, so that $T_{0}\neq 0$ and
$t_{0}T_{0}=(1-t_{0})T_{0}$,
which implies that $t_{0}=1/2$, that is not possible since $|t_{0}-1/2|>\delta$.
Then \eqref{tfinal} is true and  we also claim that there exists
a positive constant $\alpha_{0}$ such that
\beq\label{hfinal}
\intr h_{w}(p)<m-\alpha_{0}\;\;  \forall \, p\in [0,L]:\,|p-1|>\sigma_{0}, 
\eeq
where $h_{v}(p)$ is defined in \eqref{def:finalh},
Indeed, by contradiction,
there exists a sequence $\alpha_{n}$ converging to zero,
and a sequence of point $p_{n}\in [0,L]$ such that
\beq\label{hineq}
|p_{n}-1|>\sigma_{0},\quad \text{and }\quad \intr h_{w}(p_{n})\geq m-\alpha_{n}.
\eeq
Up to a subsequence, there exists $p_{0}\in [0,L]$ such that
$p_{n}\to p_{0}$. Taking limit in \eqref{hineq} we find
\beq\label{hineq2}
|p_{0}-1|>\sigma_{0},\quad \text{and }\quad \intr h_{w}(p_{0})\geq m.
\eeq
Since $\ov{p}=1$ is the unique  maximum point of $h_{w}(p)$ 
there results 
$$
m=\intr h_{w}(1) \geq \intr h_{w}(p_{0})\geq m.
$$
yielding  $p_{0}=\ov{p}=1$ contradicting
the first inequality in \eqref{hineq2}, so that \eqref{hfinal} is true.
Assume that the first inequality in \eqref{tfinal} is true, 
then we can make the choice $(p,w)=(tT^{R}_{t},w^{Ry})$
in \eqref{hfinal} yielding 
$$
I(tT^{R}_{t}w^{Ry})\leq m-\alpha_{0},\qquad
I((1-t)T^{R}_{t}w^{Rx_{0}})\leq m.
$$ 
Then, using these informations in \eqref{ineq7}, one gets
$$
\io(T^{R}_{t}U^{R}_{t})\leq 2m-\alpha_{0}+O(\eps_{R}),
$$
and an analogous argument gives the same conclusion when the second inequality in \eqref{tfinal} holds true.
This, Lemma \ref{equal}, and \eqref{ineq8} yield \eqref{eq:dim}
giving  the conclusion.

\edim

{\bf Proof of Theorems \ref{rho} and \ref{rhointro}}
\\
 Lemma \ref{geometry} yields 
 all the required  geometrical properties to apply a Linking Theorem.
Moreover from Lemmas \ref{bbound} and \ref{prop:chiave} we deduce that
\begin{align}\label{critlev}
m_{\Omega}\underset{(a)}{<}d &\underset{(b)}{<} 2m_{\Omega}.
\end{align}
Indeed, $(b)$ in \eqref{critlev} follows directly from Lemma
\ref{prop:chiave} as
$h:=id\in  {\mathcal H}$,  so that
$$
\inf_{{\mathcal H}}\sup_{Q}I_{\Omega}(h(u))\leq \sup_{Q}I_{\Omega}(u)<2m_{\Omega}.
$$
In order to show $(a)$ in \eqref{critlev}, observe that 
exploiting property \eqref{hmap} we can say that
for every  $h\in{\mathcal H}$
there exists $v\in Q$ such that $h(v)\in S$, so that
$$
\max_{u\in Q}I_{\Omega}(h(u))\geq I_{\Omega}(h(v))\geq \inf_{S}I_{\Omega}>m_{\Omega},
$$
where the last inequality is proved in Lemma \ref{bbound}.
Corollary \ref{ceramirange} implies that the (PS) condition
is satisfied at level $d$. 
Then we can apply the Linking Theorem (see e.g.  
Theorem 8.22 in \cite{amma} or Theorem 8.4  in \cite{struwe}  joint with Lemma 5.14 and 5.15 in \cite{wi}) to obtain the existence
of a constrained critical point of $I$ in $\nhlo$.
This, taking into  account Remark \ref{natu}, gives the conclusion.
\\
Theorem \ref{rhointro} follows directly from Theorem \ref{rho}
as $f$ given in \eqref{fmod} satisfies all the requested hypotheses.
\edim

\bos
The previous topological argument shows that the set $\Psi_{R}(\partial{B_R(0)})$ is contractible in the sub-level set  $I_{\Omega}^{\ov{c}}$ with $\ov{c}$ given by 
$$
\ov{c}=\sup_{\ov{B_R(0)}}I(\Psi_{R}).
$$ 
Because the whole set
$\Psi_{R}(\overline{B_R(0)})$ is contractible in  $I_{\Omega}^{\ov{c}}$, but  $\Psi_{R}(\partial{B_R(0)})$ 
is not contractible in the  sub-level set 
$I_{\Omega}^{\underline{c}}$ for $\underline{c}$ given by
$$
\underline{c}=\max_{\Psi_{R}(\partial{B_R(0)})}I.
$$
Moreover, we believe that using appropriate topological tools one can show 
also multiplicity results depending on the topology of $\Omega$.
We leave this topic as an interesting, in our opinion, open problem.
\eos

{\bf Acknowledgement.} Part of this work has been done while the second author was visiting the University of Brasilia. 
She wishes to thank all the departamento de Matem\'atica
for the warm hospitality.
A special thank goes to  the Liliane's family for the hospitality and
the friendly atmosphere.
\\
The authors thank the anonymous referee whose important comments
helped them to improve their work.

\end{document}